\documentclass[a4paper,12pt]{article}
 \usepackage[all]{xy}
 \usepackage{amsfonts}
 \usepackage{amsmath}
  \usepackage{amsthm}

 \newtheorem{mth}{Theorem}[section]
 \newtheorem{mlem}[mth]{Lemma}
 
 \newtheorem{mcor}[mth]{Corollary}
 
 \theoremstyle{definition} 
 \newtheorem{mdef}{Definition}[section]
\newtheorem{mex}[mth]{Example}

 \newcommand{\mK}{\mathbb{K}}
 \newcommand{\mR}{\mathbb{R}}
 \newcommand{\mC}{\mathbb{C}}
 
 \newcommand{\mN}{\mathbb{N}}
 \newcommand{\Mod}{\mathcal{M}}
 \newcommand{\mul}{m}
 \newcommand{\comul}{\Delta}
 \newcommand{\counit}{\varepsilon}
 
 \newcommand{\act}{\rho}
 \newcommand{\coact}{\rho}
 \newcommand{\sw}[1]{{}_{(#1)}}
 \newcommand{\co}{\text{co}}
 \newcommand{\ot}{\otimes}
 \newcommand{\can}{\text{\rm can}}
 \newcommand{\trans}{\tau}
 \newcommand{\twa}{{}^{[1]}}
 \newcommand{\twb}{{}^{[2]}}

\newcommand{\lwa}{{}^{\tilde{[1]}}}
 \newcommand{\lwb}{{}^{\tilde{[2]}}}
 \newcommand{\lcan}{\tilde{\can}}
 \newcommand{\ltrans}{\tilde{\trans}}
 \newcommand{\stl}{\ell}
 \newcommand{\swa}{{}^{\underline{[1]}}}
\newcommand{\swb}{{}^{\underline{[2]}}}
\newcommand{\eqmap}[3]{\ar@<1ex>[#1]^-{#2}  \ar@<-1ex>[#1]_-{#3}}
\newcommand{\twosid}[3]{\ar@<1ex>@{<-}[#1]^{#2}  \ar@<-1ex>[#1]_{#3}}

\begin{document}
 \title{Cotensor products of quantum principal bundles.}
 \author{Bartosz Zieli\'nski\footnote{
 Department of Mathematics, University of Wales Swansea, Singleton
Park, Swansea SA2 8PP, U.K.;
mabpz@swan.ac.uk.}{ }\footnote{
Department of Theoretical Physics II, University of \L{}\'od\'z 
Pomorska 149/153, 90-236 \L{}\'od\'z, Poland.}}
 \maketitle
 
 \begin{abstract}
 \noindent A cotensor product $A\Box_H P$ of an $H$-Hopf Galois extension $A$ and a $C$-coalgebra
 Galois extension $P$, such that $P$ is an $(H,C)$-bicomodule, is analyzed. Conditions are stated,
 when $A\Box_H P$ is a $C$-coalgebra Galois extension and when there exists a strong connection
 on $A\Box_H P$. Two examples are given, in both, $A$ and $P$ are Matsumoto spheres, and
 $H=C=C(U(1))$.
 \end{abstract}
 
 \section{Introduction}

 The recent years have observed emergence of many examples of  quantum spaces and in particular of quantum 
principal bundles (\cite{Dab:Ga} collects  many examples). One of the most important methods of 
researching mathematical structures is to define methods of constructing new spaces out of existing ones, and looking 
how the properties of the former depend on the properties of the latter. For instance, in the algebraic geometry,
given two topological spaces we can construct their smash product, wedge product, etc.

This paper was inspired by \cite{Schn:homotop}, where it was shown how the passing from a 
given quantum principal bundle 
to its prolongation preserves the homotopy class of the bundle. Cotensoring of quantum principal bundles 
is a natural generalization of prolongating of a principal bundle by a Hopf algebra (see below for details).
It may prove interesting, in future work, to check the properties of the cotensor product of quantum principal bundles
(such as homotopy classes, associated bundles, differential calculus, etc.) against the respective 
properties of its building blocks.
 
 The general idea of the method is as follows. Let $C$ be a coalgebra and $H$ a Hopf algebra.
 Assume that $A$ is an $H$-Hopf Galois extension of $B$, and $P$ is a $C$-coalgebra Galois extension of $D$ and also an $(H,C)$-bicomodule and an $H$-comodule algebra. One can form a cotensor product
 $A\Box_H P$ of $A$ and $P$. It is easy to see that $A\Box_H P$ is an algebra and a right
 $C$-comodule. It turns out that, under a number of not very restrictive conditions,
 $A\Box_H P$ is a $C$-coalgebra Galois extension of $A\Box_H D$.
 
 Classically, the cotensor product of two algebras of functions on compact Hausdorf topological spaces corresponds to the algebra of functions on the Cartesian product of spaces modulo group action.
 Explicitly, let $A=\vartheta(X)$, $P=\vartheta(Y)$, $H=\vartheta(G)$, where $G$ is a topological group, which acts on the right 
 on a topological space $X$, and on the left on a topological space $Y$, so that
 $H$ is a Hopf algebra, $A$ is a right $H$-comodule algebra and $P$ is a left $H$-comodule algebra.
 Then the group $G$ acts on the left on the  Cartesian product
 $X\times Y$:
 \[
 (X\times Y)\times G\rightarrow X\times Y,\ \  ((x,y),g)\mapsto (xg,g^{-1}y),
 \]
 and $A\Box_H P=\vartheta((X\times Y)/G)$.
 
 Suppose, that the action of $G$ on $X$ is free (i.e., that $X$ is a principal bundle). Let
 $M=X/G$ be the base space of this bundle, and let $\pi:X\rightarrow M$ be a natural projection. 
 Then $X$ is a locally trivial fibre bundle, i.e., there exists an open cover $(U_\alpha)_{\alpha\in I}$ of 
 $M$, such that $\pi^{-1}(U_\alpha)$ is homeomorphic to $U_\alpha\times G$.
 
 On the other hand it is easy to see that surjection 
 \[
 \tilde{\pi}:(X\times Y)/G\rightarrow M,\ \ (x,y)\mapsto \pi(x),
 \]
 is well defined, and, $\tilde{\pi}^{-1}(U_\alpha)=U_\alpha\times Y$. Hence $(X\times Y)/G$ is a locally trivial fibre bundle, with a base space $M$, and a fibre $Y$.
 
 Suppose that another group $K$ acts on $Y$ on the right, and actions of $G$ and $K$ on $Y$ commute (i.e., dually, for $C=\vartheta(K)$, $P$ is an $(H,C)$-bicomodule). Assume also that the action of $K$ on $Y$ is free, i.e., $Y$ is also a principal bundle, with a 
 structure group $K$. Under certain conditions
 natural action of $K$ on $(X\times Y)/G$:
 \[
 ((x,y),k)\mapsto (x,yk)
 \]
 is also free, hence $(X\times Y)/G$ is a principal bundle with structure group $K$ and base space
 $(X\times (Y/K))/G$. Then we can view $(X\times Y)/G$ in two ways as a tower of fibre bundles.
 Firstly, it is a fibre bundle with the base space $M$ and the fibre $Y$, which, in turn, is is also a fibre bundle with the base $Y/K$ and the fibre $K$.
Secondly, it is a fibre bundle with the fibre $K$ and the base $(X\times (Y/K))/G$, which in turn is also a fibre
bundle with the base $M$ and the fibre $Y/K$.

The construction of the cotensor product of quantum principal bundles can be thought of as a 
natural generalization
of prolongations of Hopf-Galois extensions (\cite{schn:Princip}).

Let $A$ be an $H$-Hopf Galois extension of $B$ and let $P$ be a Hopf algebra. Suppose that
$f:P\rightarrow H$ is a surjective Hopf algebra morphism. Define left $H$-coaction on $P$ by
$p\mapsto f(p\sw{1})\ot p\sw{2}$, where we use Sweedler notation for coproduct, 
the summation is implicitly understood.
 Then $A\Box_H P$ is a $P$-Hopf Galois extension of $B$ called
a $P$-prolongation of $A$. As $P$ is a $P$-Hopf Galois extension of the ground ring, this is a special case of the cotensor product of quantum principal bundles.

Another special case of the construction described in this paper is a cotensor product of bigalois objects
(c.f. \cite{Sch:Tft}).

The construction of the cotensor product of quantum principal bundles, described in the paper, can be a rich source of new examples of $C$-coalgebra Galois extensions, and it also offers an insight into the structure of coalgebra Galois extensions.

Organization of the paper is as follows:

{\bf Section 2} recalls basic basic terminology and notation related to 
 $C$-coalgebra Galois extensions  and  entwined modules.

{\bf Section 3} contains results concerning the existence of the inverse of the canonical map for the cotensor product
$A\Box_H P$ of an $H$-Hopf Galois extension $A$ of $B$ and a $C$-coalgebra Galois extension
$P$ of $D$ (where $P$ is an $(H,C)$-bicomodule), in the case when $D=P^{\co C}\subseteq {}^{\co H}P$, i.e., when the algebra
of functions on the base of the 
fibration $A\Box_H P$ is $B\ot D$. At the end of the section we consider  the cotensor product 
of two copies of noncommutative Matsumoto sphere, as an example.

{\bf Section 4} is concerned with the existence of an entwining and a  strong connection on the cotensor product of 
 quantum principal bundles. At the end of the section we give as an example, the cotensor product of 
two copies of the Matsumoto sphere, with the coaction defined differently than in the example in  Section 3.

\section{Basic terminology and notation}

Unless otherwise stated, we work over a general, commutative and unital ring $\mK$.

\paragraph{Categories of modules and comodules.}
Let $C$ and $H$ be coalgebras, and $A$ and $B$ algebras. We denote by ${}_A\Mod$,
$\Mod_B$, ${}^H\Mod$, $\Mod^C$, ${}_A\Mod_B$, ${}^H\Mod^C$, ${}_A\Mod^C$, etc, respectively, 
the category of left $A$-modules, right $B$-modules, left $H$-comodules, right $C$-comodules,
$(A,B)$-bimodules, $(H,C)$-bicomodules,
left $A$-modules and right $C$-comodules such that
the $C$-coaction  commutes with the $A$-action, etc.

\paragraph{Comultiplication, coaction and the Sweedler notation.} 
Suppose $C$ and $H$ are coalgebras, $M\in {}^H\Mod$ and $N\in\Mod^C$.
We denote the comultiplication by $\comul:C\rightarrow C\ot C$, the left $H$-coaction
${}^H\coact:M\rightarrow H\ot M$, the right $C$-coaction by 
$\coact^C:N\rightarrow N\ot C$.
We also use the Sweedler notation: $\comul(c)=c\sw{1}\ot c\sw{2}$, ${}^H\coact(m)=m\sw{-1}\ot m\sw{0}$,
$\coact^C(n)=n\sw{0}\ot n\sw{1}$ for all $c\in C$, $m\in M$, $n\in N$, the summation is implicitly understood.

\paragraph{Entwining structures.}  Entwining structures were introduced in \cite{BrzMaj:CoalBun}.
A good introduction to entwining structures can be found in
\cite{BrzWis:Cor}, as well as in \cite{BrzMaj:Fac} and \cite{Brz:ModAs}. Let $C$ be a coalgebra and 
$P$ an algebra. Throughout the paper, by entwining map we mean a right-right entwining map
$\psi:C\ot P\rightarrow P\ot C$, corresponding entwining structure is denoted by $(P,C)_\psi$.
We use the following summation notation for entwining map $\psi$, and, if entwining map is bijective, its inverse 
$\psi^{-1}:P\ot C\rightarrow C\ot P$:
\[
\psi(c \ot p)=p_\alpha\ot c^\alpha, \  \psi^{-1}(p\ot c)=c_A\ot p^A, \ \text{for all }c\in C,\ p\in P, 
\]
where, respectively, small Greek and large Latin letters are used for implicit summation indices.
The category of entwined modules associated to $(P,C)_\psi$ is denoted $\Mod^C_P(\psi)$.

Let $(P,C)_\psi$ be an entwining structure, and let $A$ be an algebra and an entwined module.
An algebra extension $B\subseteq A$ is called a $(P,C)_\psi$-extension iff
$B=A^{\co C}$. Of particular interests are $(P,C)_\psi$-extensions $B\subseteq P$. Such extensions
are denoted by $P(B,C,\psi)$. In this case, if there exists a grouplike element $e\in C$
such that,  for all $p\in P$, $\coact^C(p)=\psi(e\ot p)$ then $P_e(B,C,\psi)$ is called an $e$-copointed
$(P,C)_\psi$-extension.

\paragraph{The canonical map and quantum principal bundles.} Suppose that $C$ is a coalgebra and $P$ is an algebra and a right $C$-comodule. To set the notation, we recall definition of a canonical map,
\[
\can^C_P:P\ot_B P\rightarrow P\ot C,\ p\ot_B p'\mapsto pp'\sw{0}\ot p'\sw{1},
\]
where $B=P^{\co C}$ is the algebra of coinvariants of the $C$-coaction. If the canonical map is 
an isomorphism in ${}_P\Mod^C$, then
$P^C(B)$ is called a $C$-coalgebra Galois extension of $B$. If, in addition, $C$ is a Hopf algebra, and $P$ is a comodule algebra, $P^C(B)$ is called a $C$-Hopf Galois extension.
We recall the definition of a translation map,
\[
\trans^C_P:C\rightarrow P\ot_B P,\ \trans^C_P(c)=(\can^C_P){}^{-1}(1_P\ot c),\ \text{for all }c\in C.
\]
We use explicit `Sweedler-like' summation notation for the translation map:
\[
\trans^C_P(c)=c\twa\ot_B c\twb,\ \text{for all }c\in C,
\]
where implicit summation is understood.
In what follows, we use frequently the following properties of the translation map (\cite{Sch:Rep}),
for all $c\in C$, $p\in P$,
\begin{gather}
 c\twa c\twb\sw{0}\ot c\twb\sw{1}=1_P\ot c, \label{transdef}\\
c\sw{1}\twa\ot_B c\sw{1}\twb\ot c\sw{2}=c\twa\ot_B c\twb\sw{0}\ot  c\twb\sw{1},\label{colin}\\
 c\twa c\twb=\counit(c),\label{trmuleps}\\
 p\sw{0}p\sw{1}\twa\ot_B p\sw{1}\twb=1_P\ot_B p,\label{trpp}\\
 c\sw{1}\twa\ot_B c\sw{1}\twb c\sw{2}\twa\ot_B c\sw{2}\twb=c\twa\ot_B 1_P\ot_B c\twb\\
 c\twa\sw{0}\ot_B c\twb\ot c\twa\sw{1}=c\sw{2}\twa\ot_B c\sw{2}\twb\ot Sc\sw{1}.
 \label{leftc}
 \end{gather}
 If $P^C(B)$ is a $C$-Hopf Galois extension, then also, for all $b\in B$, $c,c'\in C$,
 \begin{gather}
 bc\twa\ot_B c\twb=c\twa\ot_B c\twb b,\label{trcomu}\\
 (cc')\twa\ot_B (cc')\twb=c\twa c'\twa\ot_B c'\twb c\twb,\label{trmul}.  
 \end{gather}

\section{The inverse of the canonical map for the cotensor product of quantum principal bundles.}

The general idea of the construction is as follows: take an $H$-Hopf Galois extension $A^H(B)$
 and a $C$-coalgebra Galois extension $P^C(D)$ such that $P$ is in addition a left $H$-comodule,
 and form the cotensor product $A\Box_H P$.

 \begin{mlem}\label{algcot}
  If $C$ is a coalgebra, flat as a $\mK$-module, $P$ is an $(H,C)$ bicomodule and $A$ is a right
 $H$-comodule, then $A\Box_H P$ is a right $C$-comodule. Moreover, if $H$ is a bialgebra and
 $A$ and $P$ are $H$-comodule algebras, then $A\Box_H P$ is an algebra.
 \end{mlem}
 \begin{proof} 
 The first statement is the standard result in the coalgebra theory (c.f. II3 \cite{BrzWis:Cor}). We include its
 proof for completeness. 
 
 There exists a natural right $C$-comodule structure on $A\ot P$, given  by the coaction on the 
 second factor,
 $\coact^C_{A\ot P}:a\ot p\mapsto a\ot \coact^C_P(p)$. Let
 \[
 \tilde{\coact}^C_{A\Box_H P}=\left.\coact^C_{A\ot P}\right|_{A\Box_H P}:A\Box_H P\rightarrow A\ot P\ot C.
 \]
 By the flatness of $C$ over $\mK$, the exactness of the defining sequence,
 \begin{equation}
 \xymatrix@1{
 0 \ar[r] & A\Box_H P\ar[r] & A\ot P \eqmap{rr}{\coact^H_A\ot P}{A\ot{}^H\coact_P} & & A\ot H\ot P
 }
 \end{equation}
 implies the exactness of the top row of the following diagram:
 \begin{equation}
 \xymatrix{
 0 \ar[r] & (A\Box_H P)\ot C\ar[r] & A\ot P\ot C \eqmap{rr}{\coact^H_A\ot P\ot C}{A\ot{}^H\coact_P\ot C} & & 
 A\ot H\ot P\ot C\\
                                                        &  &A\Box_H P \ar[u]_{\tilde{\coact}^C_{A\Box_H P}}
                                                        \ar@{-->}[lu]^{\coact^C_{A\Box_H P}}& &
 }.
 \label{diaco}
 \end{equation}
 As $P$ is an $(H,C)$ bicomodule, 
 \begin{multline*}
 (\coact^H_A\ot P\ot C-A\ot{}^H\coact_P\ot C)\circ\tilde{\coact}^C_{A\Box_H P}\\=
 (A\ot H\ot\coact^C_P)\circ\left.(\coact^H_A\ot P-A\ot{}^H\coact_P)\right|_{A\Box_H P}=0,
 \end{multline*}
  by the definition of
 $A\Box_H P$. Hence, by the universal property of the kernel, there is a unique factorisation 
 $\coact^C_{A\Box_H P}$ (\ref{diaco}).
 
 One defines the multiplication on $A\Box_H P$ first by restricting the usual tensor product multiplication 
 $(a\ot p)\ot(a'\ot p')\mapsto (aa'\ot pp')$ to 
 $\tilde{\mul}_{A\Box_H P}:(A\Box_H P)\ot (A\Box_H P)\rightarrow (A\ot P)$. Then, using that $A$ and $P$ are $H$-comodule algebras, one proves that it gives zero, when composed with the equalising map, so there exists a unique factorisation 
 $\mul_{A\Box_H P}:(A\Box_H P)\ot (A\Box_H P)\rightarrow (A\Box_H P)$.
 \end{proof}

The following lemma is probably well-known in the ring and module theory. However, since we were not able to find an exact reference, we carefully provide an explicit proof.

\begin{mlem}
 \label{tenfl}
 Suppose that $A$ and $B$ are algebras, $M$  is a right $A$ module and $N$ is a right $B$ module. Then $M\ot N$ is a right $A\ot B$ module, with the obvious tensor action   
 $(m\ot n)\cdot (b\ot d)\mapsto mb\ot nd$. On the other hand, if $P$ is a left $A\ot B$ module, it is also
 a left $A$ and $B$ module, with the left $A$-action $a\cdot p\mapsto (a\ot 1_B)p$, and the left $B$-
 action $b\cdot p\mapsto (1_P\ot b)p$. Obviously, these $A$ and $B$-actions commute. Furthermore 
 $(N\ot_B P)\in {}_A\Mod$, with the left $A$-action
 $a\cdot (n\ot_B p)\mapsto n\ot_B(ap)$.
 The following statements are true:
 \begin{enumerate}
 \item For all $P$ in ${}_{A\ot B}\Mod$, the map
 \begin{gather}
 \phi_P:M\ot_A(N\ot_B P)\rightarrow (M\ot N)\ot_{A\ot B} P,\nonumber\\
 m\ot_A(n\ot_B p)\mapsto (m\ot n)\ot_{A\ot B} p
 \end{gather}
 is a $\mK$-linear isomorphism, with the inverse given explicitly by:
 \begin{gather}
 \phi^{-1}_P:(M\ot N)\ot_{A\ot B} P\rightarrow M\ot_A(N\ot_B P),\nonumber\\
 (m\ot n)\ot_{A\ot B} p\mapsto m\ot_A(n\ot_B p).
 \end{gather}
 \item If $M_A$ and $N_B$ are flat,  then $M\ot N$ is a flat
right $A\ot B$-module.
\item If $M$ is a faithfully flat right $A$-module and $N$ is a 
faithfully flat right $B$-module, then $M\ot N$ is 
a faithfully flat right $A\ot B$-module.
\item If $M$ is a faithfully flat right $A$-module and $M\ot N$ is a flat right $A\ot B$ module, then $N$ is a flat right $B$-module.
\item If $M$ is a faithfully flat right $A$-module and $M\ot N$ is a faithfully flat
$A\ot B$-module, then $N$ is a faithfully flat $B$-module.
 \end{enumerate}
 \end{mlem}
 \begin{proof}
 \newcommand{\phio}{\overline{\phi_P}}
 \newcommand{\phioo}{\overline{\overline{\phi_P}}}
 \newcommand{\iphi}{\phi_P^{-1}}
 \newcommand{\iphio}{\overline{\phi_P^{-1}}}
 
 {\noindent\bf 1.} First we show that $\phi_P$ is well defined. Indeed, consider the map
 \begin{gather*}
 \phioo:M\ot(N\ot P)\rightarrow (M\ot N)\ot_{A\ot B} P,\\
 m\ot (n\ot p)\mapsto (m\ot n)\ot_{A\ot B} p.
  \end{gather*}
 The map $\phioo$ is the composition of, the associativity isomorphism from the tensor product over $\mK$ with the canonical surjection from the tensor product over $\mK$ to the tensor product over 
 $A\ot B$, hence it is well defined.
 
  Since the tensor product functor is right exact, 
  the top row in the following diagram, obtained by applying
  functor $M\ot \cdot$ to the exact sequence defining tensor product over $B$, is exact:
  \begin{equation}
  \xymatrix{
    M\ot (N\ot B \ot P)\eqmap{r}{M\ot\act_B\ot P}{M\ot N\ot {}_B\act} &  M\ot (N\ot P)\ar[r] 
    \ar[d]_{\phioo}& 
    M\ot (N\ot_B P)\ar[r]  \ar@{-->}[dl]^{\phio}& 0.\\
      &  (M\ot N)\ot_{A\ot B} P & &
  }
  \end{equation} 
  By the universal property of a cokernel, in order to prove that, $\phioo$ factorises (uniquely) through 
  $\phio$ and the canonical surjection onto the tensor product over $B$, it is enough to show that 
  \[
  \phioo\circ (M\ot\act_B\ot P)=\phioo\circ(M\ot N\ot {}_B\act). 
  \]
  Indeed, for all $m\ot (n\ot b\ot p)\in M\ot (N\ot B \ot P)$,
  \begin{multline*}
  \phioo\circ (M\ot\act_B\ot P)(m\ot (n\ot b\ot p))=\phioo(m\ot(nb\ot p))\\
  =(m\ot nb)\ot_{A\ot B} p=(m\ot n)(1_A\ot b)\ot_{A\ot B} p\\
  =(m\ot n)\ot_{A\ot B}(1_A\ot b) p=(m\ot n)\ot_{A\ot B} bp\\
  =\phioo(m\ot(n\ot bp))
  =\phioo\circ(M\ot N\ot {}_B\act)(m\ot(n\ot b\ot p)). 
  \end{multline*}
  Similarly, (in the diagram below, with $Q=N\ot_B P$),
  \begin{equation*}
  \xymatrix{
  M\ot A\ot Q\eqmap{r}{\act_A\ot Q}{M\ot{}_A\act} & M\ot Q \ar[r] \ar[d]_{\phio}
   & M\ot_A Q \ar[r] \ar@{-->}[dl]^{\phi_P} & 0\\
   & (M\ot N)\ot_{A\ot B} P, &&
  }
  \end{equation*}
  one shows that, for all $m\ot a \ot (n\ot_B p)\in M\ot A\ot (N\ot_B P)$,
  \begin{multline*}
  \phio\circ(\act_A\ot(N\ot_B P))(m\ot a \ot (n\ot_B p))=
  \phio(ma\ot (n\ot_B p))\\
  =(ma\ot n)\ot_{A\ot B} p=(m\ot n)(a\ot 1_B)\ot_{A\ot B} p=(m\ot n)\ot_{A\ot B} ap\\
  =\phio(m\ot (n\ot ap))=\phio\circ(M\ot {}_A\act)(m\ot a \ot (n\ot_B p)),
  \end{multline*}
  hence map $\phi_P$ exists, and it is a unique factorization of $\phio$, through the canonical
   surjection onto
  the tensor product over $B$.
  
  In the same way one proves that $\iphi$ is well defined. Indeed, the map
  \begin{gather}
  \iphio:(M\ot N)\ot P\rightarrow M\ot_A(N\ot_B P),\nonumber\\
  (m\ot n)\ot p\mapsto m\ot_A(n\ot_B p)
  \end{gather}
  is a composition of the associativity bijection for the  tensor product over $\mK$ 
  with two surjections onto the  tensor product over $A$ and $B$.
   Furthermore, for all $(m\ot n)\ot(a\ot b)\ot p\in(M\ot N)\ot (A\ot B)\ot P$,
  \begin{multline*}
  \iphio\circ(\act_{A\ot B}\ot P)((m\ot n)\ot(a\ot b)\ot p)=\iphio((ma\ot nb)\ot p)\\
  =ma\ot_A (nb\ot_B p)=m\ot_A(n\ot_B abp)=m\ot_A (n\ot_B ((a\ot b)p))\\
  =\iphio((m\ot n)\ot ((a\ot b)p))=\iphio\circ((M\ot N)\ot {}_{A\ot B}\act)
  ((m\ot n)\ot(a\ot b)\ot p),
  \end{multline*}
  hence map $\iphi$ exists and it is a unique factorization of $\iphio$ through the surjection 
  onto the tensor product over $A$ and $B$.
  
  Since, for all $m\ot_A(n\ot_B p)\in M\ot_A(N\ot_B P)$, $(m'\ot n')\ot_{A\ot B} p'\in (M\ot N)\ot_{A\ot B} P$,
  \begin{gather*}
  \iphi\circ\phi_P(m\ot_A(n\ot_B p))=\iphi((m\ot n)\ot_{A\ot B} p)=m\ot_A(n\ot_B p),\\
  \phi_P\circ\iphi((m'\ot n')\ot_{A\ot B} p')=\phi_P(m'\ot_A(n'\ot_B p'))=(m'\ot n')\ot_{A\ot B} p',
  \end{gather*}
 the map $\phi_P$ is a linear bijection and $\iphi$ is its inverse.
  
  Let $P,P',P''\in{}_{A\ot B}\Mod$. The following diagram will be used in the remaining part of proof:
  \begin{equation}
  \xymatrix@C=15pt{
  0 \ar[r] & M\ot_A(N\ot_B P) \ar[r] \ar[d]^{\phi_P}&
   M\ot_A(N\ot_B P') \ar[r] \ar[d]^{\phi_{P'}}&
    M\ot_A(N\ot_B P'') \ar[r] \ar[d]^{\phi_{P''}}& 0\\
    0\ar[r] & (M\ot N)\ot_{A\ot B} P \ar[r] &
    (M\ot N)\ot_{A\ot B} P' \ar[r] &
    (M\ot N)\ot_{A\ot B} P'' \ar[r] & 0.
  }\label{tencomfl}
  \end{equation}

   {\noindent\bf 2.} Suppose we are given an exact sequence of left $A\ot B$-modules:
   \begin{equation}
   \xymatrix{
   0\ar[r] & P\ar[r]^f & P'\ar[r]^g & P''\ar[r] & 0
   }
   \label{seqab}
   \end{equation}
   If the maps in the top row of diagram (\ref{tencomfl}) are $M\ot_A(N\ot_B f)$
   and $M\ot_A(N\ot_B g)$, then, since $M$ and $N$ are flat moduls, the top row sequence
   in (\ref{tencomfl}) is exact. If (\ref{tencomfl}) is a commutative diagram
   then, as verical maps are linear bijections,
   the bottom row is also exact. And the bottom horizontal maps are
   \begin{gather}
   \phi_{P'}\circ (M\ot_A(N\ot_B f))\circ\phi^{-1}_{P}=(M\ot N)\ot_{A\ot B} f,\nonumber\\
   \phi_{P''}\circ (M\ot_A(N\ot_B g))\circ\phi^{-1}_{P'}=(M\ot N)\ot_{A\ot B} g.\label{botmap} 
   \end{gather}  
   Hence $M\ot N$ is a flat $A\ot B$-module.
   
   {\noindent\bf 3.} By the previous part of the lemma, $M\ot N$ is a flat $A\ot B$ module. 
   Suppose that, given a, not necessarily exact, sequence of maps (\ref{seqab}), 
   the induced sequence in the bottom row of diagram (\ref{tencomfl}) is exact.
   Let the maps in the top row be
    \begin{gather*}
   \phi_{P'}^{-1}\circ(M\ot N)\ot_{A\ot B} f \circ\phi_{P}=(M\ot_A(N\ot_B f)),\\
   \phi_{P''}^{-1}\circ(M\ot N)\ot_{A\ot B} g\circ\phi_{P'}= (M\ot_A(N\ot_B g)),
   \end{gather*}    
    so that diagram
   (\ref{tencomfl}) is commutative. Since vertical maps are bijective and bottom row is exact, the top row will be exact as well. Then, using the faithful flatness of $M$ and $N$, one can can reconstruct the sequence (\ref{seqab}) from the top row of (\ref{tencomfl}), knowing it has to be exact.
   Hence $M\ot N$ is faithfully flat over $A\ot B$.
   
{\noindent\bf 4.} Given an exact sequence of left $A\ot B$-modules (\ref{seqab}), the 
flatness of $(M\ot N)$ implies the 
exactness of the bottom row of (\ref{tencomfl}) with injection $(M\ot N)\ot_{A\ot B} f$ and surjection
$(M\ot N)\ot_{A\ot B} g$.
 By the commutativity of (\ref{tencomfl}) and the bijectivity of vertical maps,
the top row of (\ref{tencomfl}) is also exact with injection $M\ot_A(N\ot_B f)$, and surjection
$M\ot_A(N\ot_B g)$. Then use the faithful flatness of $M$ to infer the exactness of the sequence
\begin{equation*}
\xymatrix{
0\ar[r] &  N\ot_B P\ar[r]^{N\ot_B f} & N\ot_B P' \ar[r]^{N\ot_B g} & N\ot_B P''\ar[r] & 0.
}
\end{equation*}

{\noindent\bf 5.} By the previous part of the lemma, $N$ is a flat $B$-module.
Given morphisms of left $A\ot B$-modules $f:P\rightarrow P'$ and
$g:P'\rightarrow P''$, assume that the sequence
\begin{equation}
\xymatrix@1{
0\ar[r] & N\ot_B P\ar[r]^{N\ot_B f} & N\ot_B {P'}\ar[r]^-{N\ot_B g} & N\ot_BP''\ar[r] & 0
}
\end{equation}
is exact. Then, by the 
flatness of $M$, the induced sequence in the top row of diagram (\ref{tencomfl}) will be 
exact. Define maps in the bottom row using formula (\ref{botmap}), so that (\ref{tencomfl}) is commutative, and the bottom sequence is exact.
Then use the faithful flatness of $M\ot N$ to prove the exactness of sequence (\ref{seqab}).
 \end{proof}
 
 The following two lemmas are generalisations of known properties of Hopf Galois extensions.
 
  \begin{mlem} \cite{Tak:Ext}\label{flbij}
 Let $C$ be a coalgebra and $P$  an algebra and a right $C$-comodule. Take any subalgebra
 $B\subseteq  P^{\co C}$ and define
 \begin{equation}
 \can^C_{P(B)}:P\ot_B P\rightarrow P\ot C,\ \ \  p\ot_B p'\mapsto pp'\sw{0}\ot p'\sw{1}.
 \end{equation}
 If $\can^C_{P(B)}$ is bijective and $P$ is right faithfully flat over $B$, then $B=P^{\co C}$ and
 $P^C(B)$ is a $C$-coalgebra Galois extension. We will often refer to $\can^C_{P(B)}$ as a 
 {\em canonical map}.
 \end{mlem}
 \begin{proof}
 This lemma can be viewed as a special case of Theorem 4.12 \cite{KaTor:Mor}, 
 but it can also be proven directly as follows.

 The canonical map $\can^C_{P(B)}$ maps 
 $P\ot_B P^{\co C}$ onto $P\coact^C(1_P)\subseteq P\ot C$. Moreover, 
 $P\coact^C(1_P)\simeq P$, where the isomorphism and its inverse are given by:
 \[
 p1\sw{0}\ot 1\sw{1}\mapsto p1\sw{0}\counit(1\sw{1})=p,\ \ \ p\mapsto p1\sw{0}\ot 1\sw{1}.
 \]
 On the other hand, there is an isomorphism $P\ot_B B\simeq P$ given by the formulae
 $p\ot_B b\mapsto pb$, $p\mapsto p\ot_B 1_P$. This leads to the  sequence of maps
 \begin{equation*}
 P\simeq P\ot_B B\stackrel{P\ot\imath}{\longrightarrow} P\ot_B P^{\co C}\simeq  P\coact^C(1_P)
 \simeq P,
 \end{equation*}
 whose composition is simply the 
 identity map $p\mapsto p$. Hence $P\ot \imath$ is a bijection and, by the 
 faithful flatness of $P$ over $B$, the map  $\imath$ must also be bijective. 
 \end{proof}
 
 \begin{mlem}(c.f. \cite{schn:Princip})\label{cotfflat}
 Let $A^H(B)$ be an $H$-Hopf Galois extension (i.e., $A\in{}_B\Mod_B^H$), 
 and let $V\in {}^H\Mod_D$. If $A$ is faithfully flat as a right $B$-module and $V$ is a flat right $D$-module, then $A\Box_H V$ is a flat right $B\ot D$-module. Moreover, if $V$ is a faithfully flat right $D$-module, then $A\Box_H V$ is a faithfully flat right $B\ot D$-module.
 \end{mlem}
 
 \begin{proof} $A\Box_H V$ is a right $B\ot D$ module by restriction of tensor product action
 $(A\ot V)\ot (B\ot D)\ni (a\ot v)\ot (b\ot d)\mapsto ab\ot vd\in A\ot V$. 
 There is a chain of isomorphisms
 \begin{gather}
 \xymatrix{
 A\ot_B(A\Box_H V)\simeq (A\ot_B A)\Box_H V \ar[rr]^-{\can^H_A\ot V}&&
  (A\ot H)\Box_H V\simeq A\ot V},\nonumber\\
  f:a\ot_B\sum_ia'{}_i\ot v_i\longmapsto\sum_iaa'{}_i\ot v_i.
 \end{gather}
  Call this composition $f$.
  The leftmost map in the above composition is the natural map 
  \[
   A\ot_B(A\Box_H V)\rightarrow (A\ot_B A)\Box_H V,\ \ 
  a\ot_B(\sum_ia'{}_i\ot v_i)\mapsto\sum_i(a\ot_B a'{}_i)\ot v_i,
   \]
  which, by (\cite{Tak:GeRed},\S 1) and the assumption that $A$ 
  is flat as $B$-module, is an isomorphism.
Observe that the middle transformation is well defined as $\can^H_A$ is a right $H$-colinear map.
Moreover, cotensoring with $V$ is a left exact functor and $\can^H_A$ is bijective so $\can^H_A\ot V$ is
bijective. 
View $A\ot_B(A\Box_H V)$ as a right $B\ot D$-module by the action on the second factor.
Than $f$ is a right $B\ot D$-linear map. Assertions of Lemma easily follow, using isomorphism $f$, and 
Lemma~\ref{tenfl}.
\end{proof}

 \begin{mlem}
 Let $C$ be a coalgebra and
 $H$ be a bialgebra. Let $P$ be an $(H,C)$-bicomodule and an $H$-comodule algebra,
 such that the $H$-coaction is an algebra map.
 \begin{enumerate}
 \item If $P^C(B)$ is a $C$-coalgebra Galois extension such that $B\subseteq {}^{\co H}P$, then, for any $c\in C$,
 \begin{equation}
 c\twa\sw{-1}c\twb\sw{-1}\ot c\twa\sw{0}\ot_B c\twb\sw{0}=1_H\ot c\twa\ot_B c\twb,
 \label{lefttrcov}
 \end{equation}
 where $c\twa\ot_B c\twb=\can^C_P{}^{-1}(1\ot c)$.
 \item If, in addition, $H$ is a Hopf algebra, then, for any $c\in C$,
 \begin{equation}
 c\twa\sw{-1}\ot c\twa\sw{0}\ot_B c\twb=Sc\twb\sw{-1}\ot c\twa\ot_B c\twb\sw{0}.
 \label{lefttrcovs}
 \end{equation}
 \end{enumerate}
 \end{mlem}
 \begin{proof}
 Note that conditions (\ref{lefttrcov}) and (\ref{lefttrcovs}) make sense.
 Indeed, view $H\ot P$ as a $(B,B)$-bimodule with $b\cdot (h\ot p) \cdot b'=h\ot bpb'$.
 Assumptions that a left $H$-coaction is an algebra map and $B\subseteq{}^{\co H}P$ together imply
 that the left $H$-coaction on $P$ is
 a $(B,B)$-bimodule map.
 
 {\noindent\bf 1.} This follows by applying $(H\ot(\can_P^C){}^{-1})\circ ({}^H\coact\ot C)$ 
 to both sides of identity (\ref{transdef}). Explicitly, take any $c\in C$ and compute
\begin{multline*}
 (H\ot(\can_P^C){}^{-1})\circ ({}^H\coact\ot C)(c\twa c\twb\sw{0}\ot c\twb\sw{1})\\
 \begin{split}
 {}&=(H\ot(\can_P^C){}^{-1})(
 c\twa\sw{0} c\twb\sw{0}\sw{-1}\ot c\twa\sw{0} c\twb\sw{0}\sw{0}\ot c\twb\sw{1})\\
 {}&=(H\ot(\can_P^C){}^{-1})
 ( c\twa\sw{-1} c\twb\sw{-1}\ot c\twa\sw{0} c\twb\sw{0}\sw{0}\ot c\twb\sw{0}\sw{1})\\ 
 {}&=c\twa\sw{-1} c\twb\sw{-1}\ot c\twa\sw{0} c\twb\sw{0}\sw{0}
 c\twb\sw{0}\sw{1}\twa\ot_B c\twb\sw{0}\sw{1}\twb\\
 {}&= c\twa\sw{-1} c\twb\sw{-1}\ot c\twa\sw{0}\ot_B c\twb\sw{0},
 \end{split}
 \end{multline*}
 where the last equality follows by (\ref{trpp}). On the other hand,
\begin{equation*}\begin{split} 
 (H\ot(\can_P^C){}^{-1})\circ ({}^H\coact\ot C)(1_P\ot c)&=(H\ot(\can_P^C){}^{-1})(1_H\ot 1_P\ot c)\\
 {}&=1_H\ot c\twa\ot_B c\twb,
 \end{split}
 \end{equation*}
 since $c\twa c\twb\sw{0}\ot c\twb\sw{1}=1_P\ot c$ by (\ref{transdef}), the required assertion 
 \ref{lefttrcov} follows.
 
 {\noindent\bf 2.} Apply $(\mul\ot P\ot P)\circ(S\ot{}^H\coact\ot P)$ to both sides of (\ref{lefttrcov}).
 \end{proof}

The following theorem is the main result of this section. 
  \begin{mth}\label{thbd} Suppose that:
 \begin{enumerate}
 \item $A^H(B)$ is an $H$-Hopf Galois extension such that $A$ is faithfully flat right $B$-module  and flat  as left $B$-module,
 \item $C$ is a coalgebra that is  flat as a $\mK$-module,
 \item $P^C(D)$ is a $C$-coalgebra Galois extension such that $P$ is a faithfully flat right $D$-module and a flat left $D$-module,
 \item $P$ is an $(H,C)$ bicomodule and an $H$-comodule algebra,
 \item $D\subseteq {}^{\co H}P$.
 \end{enumerate}
 Then $(A\Box_H P)^C(B\ot D)$ is a $C$-coalgebra Galois extension, and $A\Box_H P$ is a  faithfully flat right 
$B\ot D$-module. Explicitly, the inverse of the canonical map is given by:
 \begin{gather}
 (\can^C_{A\Box_H P}){}^{-1}:(A\Box_H P)\ot C\rightarrow(A\Box_H P)\ot_{B\ot D}
 (A\Box_H P),\nonumber\\
 (\sum_ia_i\ot p_i)\ot c\mapsto \sum_i(a_ic\twb\sw{-1}\twa\ot p_ic\twa)\ot_{B\ot D}
 (c\twb\sw{-1}\twb\ot c\twb\sw{0}).\label{caninv}
 \end{gather}
 \end{mth}
 \begin{proof} $A^H(B)$ and $P^C(D)$ satisfy assumptions of Lemma \ref{algcot}, hence $A\Box_H P$
 is an algebra and a right $C$-comodule. Obviously $B\ot D\subseteq(A\Box_H P)^{\co C}$, and, by
 Lemma \ref{cotfflat}, $A\Box_H P$ is  a  faithfully flat right  $B\ot D$-module, hence if $\can^C_{A\Box_H P}$ is bijective, then, by Lemma \ref{flbij}, $(A\Box_H P)^C(B\ot D)$ is a $C$-coalgebra Galois extension.
 \newcommand{\cans}{\can^C_{A\Box_H P}}
 \newcommand{\canc}{(\can^C_{A\Box_H P}){}^{-1}}
  \newcommand{\cancl}{\overline{(\can^C_{A\Box_H P}){}^{-1}}}
  \newcommand{\cancll}{\overline{\overline{(\can^C_{A\Box_H P}){}^{-1}}}}

To verify the explicit form of the inverse of the canonical map, we
 first prove that $\canc$ is well defined. Observe that the fact that $D$ is a subalgebra of ${}^{\co H}P$,
  together with assumption that $P$ is an $H$-comodule algebra, mean  that $P\in{}^H_D\Mod$, hence formula (\ref{caninv}) makes sense. Denote by 
 \[
 \cancll:(A\Box_H P)\ot C\rightarrow (A\ot P)\ot_{B\ot D}(A\ot P)
 \]
 the colifting of $\canc$, and, (for brevity) let
 $R=B\ot D$, $Q=A\ot P$. 
 Since $A$ is flat left $B$-module and $P$ is a flat left $D$-module,
 $Q$ is a flat left $R$-module by Lemma~\ref{tenfl}.
  Hence the top row in the following sequence is exact:
 \begin{equation}
 \xymatrix{
 0 \ar[r] & (A\Box_H P)\ot_R Q \ar[r]& (A\ot P)\ot_R Q 
 \eqmap{rr}{(\coact^H\ot P)\ot Q}{(A\ot {}^H\coact)\ot Q} &      & 
 (A\ot H \ot P)\ot_R Q\\
 &&&&\\
 && (A\Box_H P)\ot C \ar[uu]_{\cancll} \ar@{-->}[luu]^\cancl
 }
\end{equation}
Hence, in order to prove that there exists a unique factorisation $\cancl$ of $\cancll$ as above, it is enough, (by the universal property of the equaliser) to prove that 
\[
((\coact^H\ot P)\ot Q)\circ\cancll=((A\ot {}^H\coact)\ot Q)\circ\cancll.
\] 
 Take any $\sum_i(a_i\ot p_i)\ot c\in(A\Box_H P)\ot C$ and write,
 \begin{multline*}
 ((A\ot {}^H\coact)\ot Q)\circ\cancll(\sum_i(a_i\ot p_i)\ot c)\\
 =\sum_i(a_i c\twb\sw{-1}\twa\ot p_i\sw{-1}c\twa\sw{-1}\ot p_i\sw{0}c\twa\sw{0})\ot_R
 (c\twb\sw{-1}\twb\ot c\twb\sw{0})\\
 \text{[ use eq. (\ref{lefttrcovs})  and that   $\sum\nolimits_ia_i\ot p_i\sw{-1}\ot p_i\sw{0}=
 \sum\nolimits_ia_i\sw{0}\ot a_i\sw{1}\ot p_i$ ] }\\
 =\sum_i(a_i\sw{0}c\twb\sw{-1}\twa\ot a_i\sw{1} Sc\swb\sw{-2}\ot p_i c\twa)\ot_R
  (c\twb\sw{-1}\twb\ot c\twb\sw{0})\\
  \text{[ use eq. (\ref{leftc}) for $c\twb\sw{-1}$ ]}\\
  =\sum_i(a_i\sw{0} c\twb\sw{-1}\twa\sw{0}\ot a_i\sw{1} c\twb\sw{-1}\twa\sw{1}\ot
  p_i c\twa)\ot_R  (c\twb\sw{-1}\twb\ot c\twb\sw{0})\\
  =((\coact^H\ot P)\ot Q)\circ\cancll((\sum_i a_i\ot p_i)\ot c).
 \end{multline*}
 
 Denote $Q'=A\Box_H P$. As $A$ is right faithfully flat over $B$, and $P$ is right faithfully flat over $D$,
 then $Q'$ is right faithfully flat, hence flat, over $R$ by Lemma~\ref{cotfflat}. Therefore the top row in the diagram below is exact:
 \begin{equation*}
 \xymatrix{
 0 \ar[r] & Q'\ot_R(A\Box_H P)\ar[r] & Q'\ot_R (A\ot P) 
 \eqmap{rr}{Q'\ot(\coact^H\ot P)}{Q'\ot(A\ot{}^H\coact)}& &
 Q'\ot_R(A\ot H\ot P).\\
 & & & & \\
 & &(A\Box_H P)\ot C \ar[uu]_\cancl  \ar@{-->}[luu]^{\canc} & &
 }
 \end{equation*}  
 Moreover,  for any $\sum_i(a_i\ot p_i)\ot c\in(A\Box_H P)\ot C$,
  \begin{multline*}
  (Q'\ot(\coact^H\ot P))\circ \cancl(\sum_i(a_i\ot p_i)\ot c)\\
  =\sum_i(a_ic\twb\sw{-1}\twa\ot p_ic\twa)\ot_R
  (c\twb\sw{-1}\twb\sw{0}\ot c\twb\sw{-1}\twb\sw{1}\ot c\twb\sw{0})\\
  \text{[ use the right $H$-colinearity of the translation map for $A^H(B)$ (eq. \ref{colin}) ]}\\
  =\sum_i(a_ic\twb\sw{-2}\twa\ot c\twa)\ot_R
  (c\twb\sw{-2}\twb\ot c\twb\sw{-1}\ot c\twb\sw{0})\\
  =(Q'\ot(A\ot{}^H\coact))\circ\cancl(\sum_i(a_i\ot p_i)\ot c).
  \end{multline*}
 By the universal property of an equaliser, there exists a unique factorization
 $\canc$ of $\cancl$ as required.
 
 It remains to  prove that $\canc$ is the inverse of $\cans$. For any
 $\sum_i(a_i\ot p_i)\ot c\in(A\Box_H P)\ot C$,
 \begin{multline*}
 \cans\circ\canc(\sum_i(a_i\ot p_i)\ot c)=\\
 \begin{split}
 {}&=\sum_ia_ic\twb\sw{-1}\twa c\twb\sw{-1}\twb\ot p_ic\twa c\twb\sw{0}\ot c\twb\sw{1}\\
 \text{[use (\ref{trmuleps})]}&=\sum_i a_i\counit_H(c\twb\sw{-1})\ot p_i c\twa c\twb\sw{0}\ot c\twb\sw{1}\\
 {}&=\sum_i a_i\ot p_i c\twa c\twb\sw{0}\ot c\twb\sw{1}
 =\sum_i (a_i\ot p_i) \ot c.
 \end{split}
 \end{multline*}
 On the other hand, for any $\sum_{ij}(a_i\ot p_i)\ot_R (a'{}_j\ot p'{}_j)\in (A\Box_H P)\ot_R(A\Box_H P)$,
 \begin{multline*}
 \canc\circ\cans(\sum_{ij}(a_i\ot p_i)\ot_R (a'{}_j\ot p'{}_j))\\
 =\sum_{ij}(a_ia'{}_jp'{}_j\sw{1}\twb\sw{-1}\twa\ot p_ip'{}_j\sw{0}p'{}_j\sw{1}\twa)\ot_R
 (p'{}_j\sw{1}\twb\sw{-1}\twb\ot p'{}_j\sw{1}\twb\sw{0})\\
 \text{[ use (\ref{trpp}) ] }=\sum_{ij}(a_ia'{}_jp'{}_j\sw{-1}\twa\ot p_i)\ot_R
 (p'{}_j\sw{-1}\twb\ot p'{}_j\sw{0})\\
 \text{[use that $\sum\nolimits_{j}a'{}_j\sw{0}\ot a'{}_j\sw{1}\ot p'{}_j=
 \sum\nolimits_ja'{}_j\ot p'{}_j\sw{-1}\ot p'{}_j\sw{0}$ ]}\\
 =\sum_{ij}(a_ia'{}_j\sw{0}a'{}_j\sw{1}\twa\ot p_i)\ot_R
 (a'{}_j\sw{1}\twb\ot  p'{}_j)\\
 \text{[ use (\ref{trpp}) ]}=\sum_{ij}(a_i\ot p_i)\ot_R (a'{}_j\ot p'{}_j).
 \end{multline*}
 Thus we conclude that $\canc$ is the inverse of $\cans$ as stated. This completes the proof of the theorem.
 \end{proof}
 
 The assumption about the faithful flatness of $A$ and $P$,
 although quite restrictive,  is a usual assumption made for Galois-type extensions to be able to view 
 them as bona fide generalisations of torsors or principal bundles. Indeed, if one wants to develop
 differential geometry on Galois-type extensions in terms of strong connections, the faithful flatness
 becomes necessary (c.f.\ Theorem 2.5 \cite{BrzHaj:Chern}).
 
\begin{mex}(\cite{schn:Princip} Remark 3.11 (2)) 
Let $A$ be an $H$-Hopf Galois extension of $B$, faithfully flat as a right $B$-module,
 and let $P$ be a Hopf algebra faithfully flat as a 
$\mK$-module. Suppose that
$f:P\rightarrow H$ is a Hopf algebra morphism. Define a left $H$-coaction on $P$ by
$p\mapsto f(p\sw{1})\ot p\sw{2}$. As $D=P^{\co P}=\mK1_P$, and $\mK1_P\subseteq {}^{\co H}P$,
all the assumptions of Theorem~\ref{thbd} are satisfied, and so $A\Box_H P$ is a $P$-Hopf Galois extension of $B\ot \mK1_P\simeq B$. This $P$-Hopf galois extension is called
a $P$-prolongation of $A$ if $f$ is surjective. 
\end{mex}

\newcommand{\uone}{C^0(U(1))}
 
 \begin{mex}
 \label{matex}
 {\em The Matsumoto sphere} (\cite{Mat:Sph}) $C^0_\theta(S^3)$ is a $\ast$-algebra
generated by $a$, $b$ and relations
\begin{equation}
aa^\ast=a^\ast a,\ \ bb^\ast=b^\ast b,\ \ ab=\lambda ba,\ \ ab^\ast=\bar{\lambda}b^\ast a,\ \
aa^\ast+bb^\ast=1,
\end{equation}
 where $\lambda=e^{2\pi i\theta}$, $\theta\in\mR$. Furthermore, denote by $C^0(U(1))$, the $\ast$-Hopf
 algebra generated by unitary and grouplike $u$, i.e., 
 \begin{equation}
 uu^\ast=u^\ast u=1,\ \ \comul(u)=u\ot u,\ \ Su=u^\ast,\ \  \counit(u)=1.
 \end{equation}
 $C^0_\theta(S^3)$ is a $\ast$-$C^0(U(1))$-comodule algebra with
 \begin{equation}
 \coact^{\uone}(a)=a\ot u,\ \ \coact^{\uone}(b)=b\ot u.
 \end{equation}
 It is easy to see that the algebra of coinvariants $B=C_\theta^0(S^3)^{\co\uone}$
 is a commutative $\ast$-algebra generated by $z=aa^\ast$,
 $x_{+}=ba^\ast$, $x_{-}=(x_{+})^\ast=ab^\ast$, with an additional relation
 \begin{equation}
 z^2+x_{+}x_{-}=1,
 \end{equation}
 and hence it can be identified with the $\ast$-algebra, $C^0(S^2)$, of polynomial functions on the 
 two-sphere.
 
 $C^0_\theta(S^3)^{\uone}(C^0(S^2))$ is a $\uone$-Hopf Galois extension with the 
 inverse of the canonical map given by, for any $p\in C^0_\theta(S^3)$ and $n\in \mN\cup \{0\}$,
 \begin{gather}
 \can^{\uone}_{C^0_\theta(S^3)}(p\ot u^n)=
 \sum_{m=0}^n {n\choose m} pb^\ast{}^ma^\ast{}^{n-m}\ot_B a^{n-m} b^m,\\
 \can^{\uone}_{C^0_\theta(S^3)}(p\ot u^\ast{}^n)=
  \sum_{m=0}^n {n \choose m}pb^ma^{n-m}\ot_B a^\ast{}^{n-m} b^\ast{}^m.
 \end{gather}

 \newcommand{\mata}{C^0_\theta(S^3)^{\uone}(C^0(S^2))}
 \newcommand{\matb}{C^0_{\theta'}(S^3)^{\uone}(C^0(S^2))}
 Consider two copies of the 
 Matsumoto sphere (though with not necessarily equal deformation parameters).
  For brevity, let $A=C^0_\theta(S^3)$,
  $P=C^0_{\theta'}(S^3)$, $H=\uone$.
  By Theorem~2.5 \cite{BrzHaj:Chern} and the existence of strong connection on the Matsumoto sphere  (\cite{BrzSit:Mat}),
   $A$ and $P$ are  faithfully flat as  right and left
  $B$-modules.
  
  Define the $\ast$-comodule algebra left $H$-coaction on $P$,
  \begin{equation}
  {}^{H}\coact(a)=u^{\ast}\ot a\ \ {}^{H}\coact(b)=u^{\ast}\ot b.
  \end{equation}
  Note that, $B={}^{\co H}P$. 
  
 The cotensor product $A\Box_H P$ is a $\ast$-algebra generated by
  \begin{gather}
  \alpha=a\ot a^\ast,\ \ \beta=b\ot b^\ast, \ \  \gamma=a\ot b^\ast,\ \ \delta=b\ot a^\ast,\nonumber\\
  \alpha^\ast=a^\ast\ot a,\ \ \beta^\ast=b^\ast\ot b, \ \  \gamma^\ast=a^\ast\ot b,\ \ \delta^\ast=b^\ast\ot a,
  \label{matmatg}
  \end{gather}
  satisfying commutation relations
  \begin{align}
  \alpha\alpha^\ast&=\alpha^\ast\alpha,&
  \beta\beta^\ast&=\beta^\ast\beta,&
  \gamma\gamma^\ast&=\gamma^\ast\gamma,&
  \delta\delta^\ast&=\delta^\ast\delta,\nonumber\\
  \alpha\beta&=\lambda\lambda'\beta\alpha,&
  \alpha\beta^\ast&=\bar{\lambda}\bar{\lambda'}\beta^\ast \alpha,&
  \alpha\gamma&=\lambda'\gamma\alpha,&
  \alpha\gamma^\ast&=\bar{\lambda'}\gamma^\ast\alpha,\nonumber\\
  \alpha\delta&=\lambda\delta\alpha,&
    \alpha\delta^\ast&=\bar{\lambda}\delta^\ast\alpha,&
    \beta\gamma&=\bar{\lambda}\gamma\beta,&
        \beta\gamma^\ast&=\lambda\gamma^\ast\beta,\nonumber\\
        \beta\delta&=\bar{\lambda'}\delta\beta,&
                \beta\delta^\ast&=\lambda'\delta^\ast\beta,&
\gamma\delta&=\lambda\bar{\lambda'}\delta\gamma,&
\gamma\delta^\ast&=\bar{\lambda}\lambda'\delta^\ast\gamma,
  \end{align}
  and, in addition,
  \begin{gather}
  \alpha^\ast\alpha+\beta^\ast\beta+\gamma^\ast\gamma+\delta^\ast\delta=1,\\
  \alpha\beta=\lambda'\gamma\delta,
  \end{gather}
  where $\lambda=e^{2\pi i\theta}$ and $\lambda'=e^{2\pi i\theta'}$. By Theorem~\ref{thbd}
  $(A\Box_H P)^H(B\ot B)$ is an $H$-Hopf Galois extension. 
  In terms of generators (\ref{matmatg}), the translation map is explicitly given by
  \begin{gather}
  \begin{split}
  \trans^H_{A\Box_H P}(u^n)&=
  \sum_{p=0}^n\sum_{m=0}^{p-1}{n\choose p}{n \choose m}
  \alpha^{n-p}\delta^{p-m}\beta^m\ot_{B\ot B}
  \beta^\ast{}^{m}\delta^\ast{}^{p-m}\alpha^\ast{}^{n-p}\\
  &\quad+\sum_{p=0}^{n}\sum_{m=p}^n
  {n\choose p}{n\choose m}\alpha^{n-m}\gamma^{m-p}\beta^p\ot_{B\ot B}
  \beta^\ast{}^p\gamma^\ast{}^{m-p}\alpha^\ast{}^{n-m},
  \end{split}\\
  \begin{split}
  \trans_{A\Box_H P}^H(u^{-n})&=
  \sum_{p=0}^n\sum_{m=0}^{p-1}{n\choose p}{n \choose m}
  \alpha^\ast{}^{n-p}\delta^\ast{}^{p-m}\beta^\ast{}^m\ot_{B\ot B}
  \beta^{m}\delta^{p-m}\alpha^{n-p}\\
  &\quad+\sum_{p=0}^{n}\sum_{m=p}^n
  {n\choose p}{n\choose m}\alpha^\ast{}^{n-m}\gamma^\ast{}^{m-p}\beta^\ast{}^p\ot_{B\ot B}
  \beta^p\gamma^{m-p}\alpha^{n-m},
  \end{split}
  \end{gather}
  for any $n\in\mN$.
  Observe that the algebra of coinvariants 
  $B\ot B$ is generated by
  \begin{align}
  z\ot 1&=\alpha^\ast\alpha+\gamma^\ast\gamma,&
  x_{+}\ot 1&=\delta\alpha^\ast+\beta\gamma^\ast,&
  x_{-}\ot 1&=\alpha\delta^\ast+\gamma\beta^\ast,\nonumber\\
  1\ot z&=\alpha^\ast\alpha+\delta^\ast\delta,&
  1\ot x_{+}&=\gamma^\ast\alpha+\beta^\ast\delta,&
  1\ot x_{-}&=\alpha^\ast\gamma+\delta^\ast\beta.
  \end{align}
 \end{mex}
 
 \section{Entwining structures and strong connections.}
 
 The previous section was concerned with the cotensor product $A\Box_H P$ of an $H$-Hopf Galois extension 
$A$ and $C$-coalgebra Galois extension $P$ such that $P$ is an $(H,C)$-bicomodule
 and $P^{\co C}\subseteq {}^{\co H}P$. In the present section we drop the latter assumption. The formula
 for he inverse of the canonical map (\ref{caninv}) becomes now badly defined, because the left coaction of $H$ on $P$
no longer commutes with the multiplication by elements of $P^{\co C}$. Hence we need a quantity associated with the inverse 
of the canonical map, but without tensor products over the subalgebra of coinvariants.
 
  \begin{mdef}\label{colifflem}
 Let $C$ be a coalgebra and $P$ be an algebra and right $C$-comodule.
 Let $\lcan^C_P:P\ot P\rightarrow P\ot C$, $p\ot p'\rightarrow pp'\sw{0}\ot p'\sw{1}$ be a lifting of a canonical map, i.e., $\lcan^C_P=\can^C_{P(\mK)}$ (c.f. Lemma~\ref{flbij}).
 A linear morphism $\ltrans^C_P:C\rightarrow P\ot P$, such that, for all $c\in C$, 
 $\lcan^C_P(\ltrans^C_P(c))=1_P\ot c$, is called 
   \emph{a colifting of  translation map}. 
    For convenience in computations we introduce the notation $\ltrans^C_P(c)=c\lwa\ot c\lwb$ (summation understood).
  \end{mdef}

A colifting of a translation map which is normalized and left and right covariant is closely related to a strong connection form (c.f.\  \cite{BrzHaj:Chern}). Explicitly, we have
   \begin{mdef}
 Let $P_e(B,C,\psi)$ be an $e$-copointed $(P,C)_\psi$-extension with a bijective entwining. A 
 \emph{strong connection form} is a linear map $\stl:C\rightarrow P\ot P$ such that
 \begin{gather}
 \stl(e)=1_P\ot1_P \label{str1},\\
 \lcan^C_P(\stl(c))=1_P\ot c \text{ for all } c\in C\label{str2},\\
 (P\ot\coact^C)\circ\stl=(\stl\ot C)\circ\comul\label{str3},\\
 ({}^{C_\psi}\coact\ot P)\circ\stl=(C\ot\stl)\circ\comul\label{str4},
 \end{gather}
 where $\coact^C:P\rightarrow P\ot C$ is a right $C$-coaction and ${}^{C_\psi}\coact:P\rightarrow C\ot P$, $p\mapsto\psi^{-1}(p\ot e)$ is a left $C$-coaction, induced by the inverse of entwining 
 (c.f. \cite{BrzMaj:Fac}).
 The action of $\stl$ on elements will be denoted by  $\stl(c)=c\swa\ot c\swb$ 
 (implicit summation understood).
 
 If $C$ is a Hopf algebra with a bijective antipode and $P$ is a right $C$-comodule algebra then, unless otherwise stated, regardless of existence of entwining we shall call a linear map 
 $\stl:C\rightarrow P\ot P$ a strong connection form if it satisfies conditions (\ref{str1}-\ref{str3})
 (with $e=1_C$) and a condition (for any $c\in C$):
 \begin{equation}
 \label{str5}
 c\swa\sw{1}\ot c\swa\sw{0}\ot c\swb=Sc\sw{1}\ot c\sw{2}\swa\ot c\sw{2}\swb
 \end{equation}
 \end{mdef}
 Observe that if $P^C(B)$ is a $C$-Hopf Galois extension, then the assumption of the bijectivity of the 
 antipode is equivalent to the bijectivity of the (unique) canonical entwining, and then, conditions
 (\ref{str4}) and (\ref{str5}) are equivalent.
 
 Recall from \cite{BrzHaj:Chern}, 
 that an $e$-copointed $C$-coalgebra Galois extension $P^C_e(B)$ is called 
 \emph{a principal extension}
 if the canonical entwining $\psi_\can$ is bijective and there exists a strong connection on $P$. 

\newcommand{\ltransc}{\ltrans^C_{A\Box_H P}}
 \newcommand{\ltransa}{\ltrans^H_{A}}
 \newcommand{\ltransp}{\ltrans^C_{P}}
 
 The following lemma is the core of the results of this section.
 
 \begin{mlem}\label{coliftlem}
 Let $H$ be a Hopf algebra and $C$ a $\mK$-flat coalgebra. Let $A$ be a right $H$-comodule algebra,
 and $P$ be a left $H$-comodule algebra such that it is an $(H,C)$-bicomodule. Assume
 that $A\ot P$ and $A\Box_H P$ are flat $\mK$-modules and that
 \begin{enumerate}
 \item there exists on $A$ a colifting of the translation map $\ltransa:H\rightarrow A\ot A$, 
 $h\mapsto h\lwa\ot h\lwb$ such that, for any $h\in H$,
 \begin{gather}
  (A\ot\coact^H)\circ\ltransa(h)=\ltransa(h\sw{1})\ot h\sw{2} 
 \text{ (right covariance)},\label{rcovcolift}\\
 h\lwa\sw{1}\ot h\lwa\sw{0}\ot h\lwb=Sh\sw{1}\ot h\sw{2}\lwa\ot h\sw{2}\lwb
 \text{ (left covariance)};\label{lcovcolift}
 \end{gather}
 \item there exists on $P$ a colifting of the translation map $\ltransp:C\rightarrow P\ot P$, 
 $c\mapsto c\lwa\ot c\lwb$ such that, for any $c\in C$,
 \begin{equation}
 c\lwa\sw{-1}\ot c\lwa\sw{0}\ot c\lwb=Sc\lwb\sw{-1}\ot c\lwa\ot c\lwb\sw{0}.\label{lhcovcol}
 \end{equation}
 \end{enumerate}
Then
 $A\Box_H P$ is an algebra and right $C$-comodule and 
 the map 
\begin{gather}
\ltransc:C\rightarrow (A\Box_H P)\ot (A\Box_H P),\nonumber\\
c\mapsto (c\lwb\sw{-1}\lwa\ot c\lwa)\ot (c\lwb\sw{-1}\lwb\ot c\lwb\sw{0})\label{coliftcot}
\end{gather}
is a colifting of the translation map on $A\Box_H P$, i.e.,
\[
\lcan^C_{A\Box_H P}\circ\ltransc(c)=(1_A\ot 1_P)\ot c
\]
 \end{mlem}
 \begin{proof} $A\Box_H P$ is an algebra and a right $C$-comodule by Lemma~\ref{algcot}. 
 Since, by assumption,
$A\ot P$ and $A\Box_H P$ are flat over $\mK$, 
in order to prove that the map $\lcan^C_{A\Box_H P}$, defined by (\ref{coliftcot}), has its image in 
$(A\Box_H P)\ot (A\Box_H P)$,
it is enough 
to show that 
\[((\coact^H\ot P)\ot (A\ot P))\circ\ltransc=((A\ot{}^H\coact)\ot(A\ot P))\circ\ltransc\] 
and
\[((A\ot P)\ot (\coact^H\ot P))\circ\ltransc=((A\ot P)\ot(A\ot{}^H\coact))\circ\ltransc,\]
(c.f. proof of the Theorem~\ref{thbd}).

Take any $c\in C$, and compute
\begin{multline*}
((\coact^H\ot P)\ot (A\ot P))\circ\ltransc(c)\\
\begin{split}
{}&=(c\lwb\sw{-1}\lwa\sw{0}\ot c\lwb\sw{-1}\lwa\sw{1}\ot c\lwa)\ot
(c\lwb\sw{-1}\lwb\ot c\lwb\sw{0})  \\
\text{[use (\ref{lcovcolift})] }&=
(c\lwb\sw{-1}\sw{2}\lwa\ot Sc\lwb\sw{-1}\sw{1}\ot c\lwa)\ot
(c\lwb\sw{-1}\sw{2}\lwb\ot c\lwb\sw{0})\\
{}&=(c\lwb\sw{0}\sw{-1}\lwa\ot Sc\lwb\sw{-1}\ot c\lwa)\ot
(c\lwb\sw{0}\sw{-1}\lwb\ot c\lwb\sw{0}\sw{0})\\
\text{[ use \ref{lhcovcol} ] }&=
(c\lwb\sw{-1}\lwa\ot c\lwa\sw{-1}\ot c\lwa\sw{0})\ot
(c\lwb\sw{-1}\lwb\ot c\lwb\sw{0})\\
{}&=((A\ot{}^H\coact)\ot(A\ot P))\circ\ltransc(c).
\end{split}
\end{multline*}
 Furthermore,
 \begin{multline*}
 ((A\ot P)\ot (\coact^H\ot P))\circ\ltransc(c)\\
 \begin{split}
 {}&=(c\lwb\sw{-1}\lwa\ot c\lwa)\ot
 (c\lwb\sw{-1}\lwb\sw{0}\ot c\lwb\sw{-1}\lwb\sw{1}\ot c\lwb\sw{0})\\
 \text{[ use (\ref{rcovcolift}) ] }&=
 (c\lwb\sw{-2}\lwa\ot c\lwa)\ot
 (c\lwb\sw{-2}\lwb\ot c\lwb\sw{-1}\ot c\lwb\sw{0})\\
 {}&=((A\ot P)\ot(A\ot{}^H\coact))\circ\ltransc(c).
 \end{split}
 \end{multline*}
 Finally we check whether (\ref{coliftcot}) is a colifting of the translation map for $A\Box_H P$.
 Note that $\mul\circ\ltransa(c)=\counit(c)$, for any $c\in C$, hence
 \begin{equation*}
 \begin{split}
 \lcan^C_{A\Box_H P}\circ\ltransc(c)
 &=c\lwb\sw{-1}\lwa c\lwb\sw{-1}\lwb\ot c\lwa c\lwb\sw{0}\sw{0}\ot c\lwb\sw{0}\sw{1}\\
 {}&=\counit(c\lwb\sw{-1})1_A\ot c\lwa c\lwb\sw{0}\sw{0}\ot c\lwb\sw{0}\sw{1}\\
 {}&=(1_A\ot 1_P)\ot c.
 \end{split}
 \end{equation*}
 This completes the proof of the lemma.
 \end{proof}

 We want to study the rules governing the existence of a strong connection on $A\Box_H P$,
 depending on the existence and properties of the strong connection forms on $A$ and $P$.
 However, for the definition of a strong connection form to make sense, we first need an entwining.

 Let $(P,C)_\psi$ be an entwining structure. Assume that $P$ is a left $H$-comodule for a coalgebra $H$. We say, that entwining $\psi$ \emph{commutes with the left $H$-coaction }
 ${}^H\coact:P\rightarrow H\ot P$, if, for any $c\in C$, the $\mK$-module  map 
 $\psi(c\ot \cdot):P\rightarrow P\ot C$ is left
 $H$-colinear. 
 
 Define left $H$-comodule structure on $C\ot P$ by $c\ot p\mapsto p\sw{-1}\ot c \ot p\sw{0}$
 Let $\psi$ be bijective. We say that the inverse of entwining map, $\psi^{-1}$, 
 \emph{commutes with the left $H$-coaction}, if, for all $c\in C$, the map
 $\psi^{-1}(\cdot\ot c):P\rightarrow C\ot P$ is left $H$-colinear. 
 Note that if the  entwining map $\psi$ is bijective and commutes with the $H$-coaction, then
 also the inverse map, $\psi^{-1}$, commutes with the $H$-coaction.
 
 \begin{mlem}\label{hcent}
  Let $H$ be a bialgebra, and $(P,C)_\psi$ be an entwining structure such that $P$ is an entwined module. Suppose that $P$ is a left $H$-comodule algebra such that
  $\psi$ commutes with the
  $H$-coaction. Then $P$ is an $(H,C)$-bicomodule if and only if
  \begin{equation}\label{unitcon}
  ({}^H\coact\ot C)\circ\coact^C(1_P)=1_H\ot\coact^C(1_P).
  \end{equation}
  \end{mlem}
  \begin{proof} 
  Since $P\in\Mod^C_P(\psi)$, the right coaction  $\coact^C:P\rightarrow P\ot C$ is necessarily given by 
  $p\mapsto 1\sw{0}\psi(1\sw{1}\ot p)$.  
  Hence, for all $p\in P$,
  \begin{multline*} 
  ({}^H\coact\ot C)\circ\coact^C(p)\\
  \begin{split}
  {}&={}^H\coact(1\sw{0}p_\alpha)\ot 1\sw{1}{}^\alpha
  &\text{(since $P\in\Mod^C_P(\psi)$)}\\
  {}&={}^H\coact(1\sw{0}){}^H\coact(p_\alpha)\ot 1\sw{1}{}^\alpha 
  &\text{(since ${}^H\coact$ is algebraic)}\\ 
  {}&=1_Hp_\alpha\sw{-1}\ot 1\sw{0}p_\alpha\sw{0}\ot 1\sw{1}{}^\alpha
 &\text{(by (\ref{unitcon}))}\\
  {}&=p\sw{-1}\ot 1\sw{0}p\sw{0}{}_\alpha\ot 1\sw{1}{}^\alpha
  &\text{(since $\psi$ commutes with $H$-coaction)}\\
  {}&=(H\ot\coact^C)\circ{}^H\coact(p) 
  &\text{(since $P\in\Mod^C_P(\psi)$)}.
  \end{split}
  \end{multline*}
  Obviously, if $P$ is an $(H,C)$-bicomodule, then (\ref{unitcon}) is satisfied,
  since ${}^H\coact(1_P)=1_H\ot 1_P$.
 \end{proof}

  \begin{mlem}\label{hcent2}
  Let $P(B,C,\psi)$ be a $C$-$(P,C)_\psi$-extension and let $H$ be a Hopf algebra. Let $P$ be an 
  $(H,C)$-bicomodule and an $H$-comodule algebra, where $H$ is a Hopf algebra. 
  Suppose that either:
  \begin{description}
  \item[(a)] there exists a colifting of the translation map $\ltrans^C_P:c\mapsto c\lwa\ot c\lwb$ such that, for any $c\in C$,
  \begin{equation}
  c\lwa\sw{-1}c\lwb\sw{-1}\ot c\lwa\sw{0}\ot c\lwb\sw{0}=1_H\ot c\lwa\ot c\lwb,\label{coliftcond}
  \end{equation}
  \end{description}
  or
  \begin{description}
  \item[(b)] $P^C(B)$ is a $C$-coalgebra Galois extension and $B\subseteq{}^{\co H}P$.
    \end{description}
    Then the entwining $\psi$ commutes with the $H$-coaction ${}^H\coact:P\rightarrow H\ot P$.
  \end{mlem}
  \begin{proof}   
  Assume that there exists a colifting satisfying condition (\ref{coliftcond}).
  Since $P$ is an $(H,C)$-bicomodule,
  \begin{equation}
  \label{entbic}
  p\sw{-1}\ot 1_P\sw{0}\psi(1_P\sw{1}\ot p\sw{0})=p_\alpha\sw{-1}\ot 1_P\sw{0}p_\alpha\sw{0}\ot 1_P\sw{1}{}^\alpha
  \end{equation}
   for any $p\in P$ and $c\in C$. 
   By the definition of a colifting, $c\lwa c\lwb\sw{0}\ot c\lwb\sw{1}=1_P\ot c$, hence we can write,
   \[
   \psi(c\ot p)=c\lwa c\lwb\sw{0}\psi(c\lwb\sw{1}\ot p)=c\lwa\coact^C(c\lwb p)=
   c\lwa 1_P\sw{0}\psi(1_P\sw{1}\ot c\lwb p),
   \]
   where the second and third equalities are consequences of the fact that $P$ is an entwined module.
    Take any $c\in C$, $p\in P$, and compute
   \begin{multline*}
   p\sw{-1}\ot\psi(c\ot p\sw{0})=p\sw{-1}\ot c\lwa 1_P\sw{0}\psi(1_P\sw{1}\ot c\lwb p\sw{0})\\
    \text{[ use (\ref{coliftcond}) ]}=c\lwa\sw{-1}c\lwb\sw{-1}p\sw{-1}\ot c\lwa\sw{0}1_P\sw{0}
    \psi(1_P\sw{1}\ot c\lwb\sw{0}p\sw{0})\\
    =c\lwa\sw{-1}(c\lwb p)\sw{-1}\ot c\lwa\sw{0}1_P\sw{0}
    \psi(1_P\sw{1}\ot (c\lwb p)\sw{0})\\
    \text{[use (\ref{entbic})] }=c\lwa\sw{-1}(c\lwb p)_\alpha\sw{-1}\ot c\lwa\sw{0}1_P\sw{0}
    (c\lwb p)_\alpha\sw{0}\ot 1_P\sw{1}{}^\alpha\\
    \text{[ ${}^H\coact$
    is algebraic,  hence $1_H\ot 1_P\sw{0}\ot 1_P\sw{1}=1_P\sw{-1}\ot 1_P\sw{0}\ot 1_P\sw{1}$ ]}\\
    =(c\lwa 1_P\sw{0} c\lwb{}_\alpha p_\beta)\sw{-1}\ot (c\lwa 1_P\sw{0} c\lwb{}_\alpha p_\beta)\sw{0} \ot 1_P\sw{1}{}^{\alpha\beta}\\
    \text{[ notice that $c\lwa 1_P\sw{0}c\lwb{}_\alpha\ot 1_P\sw{1}{}^\alpha=c\lwa\coact^C(c\lwb)=1\ot c$ ]}\\
     = p_\beta\sw{-1}\ot p_\beta\sw{0}\ot c^\beta.
   \end{multline*}
   Note the similarity of the condition (\ref{coliftcond}) to the property (\ref{lefttrcov}) of the translation map.
   Indeed, (\ref{lefttrcov}) has the same form as (\ref{coliftcond}), but with adorned tensor products.
   Thus a similar argument, to the one proving the first part of the lemma, can be used to show, that hypothesis {\bf (b)} implies the assertion.
   \end{proof}
   
   Conditions (\ref{coliftcond}) or $B\subseteq{}^{\co H}P$ are not necessary. 
   For instance if $C$ is a Hopf algebra and $P^C(B)$ is a
   $C$-Hopf Galois extension, then, for the canonical entwining associated to $P^C(B)$,
   $\psi_{\can}(c\ot p)=\can^C_P((\can^C_P)^{-1}(1\ot c)p)=p\sw{0}\ot p\sw{1}c$  to commute with right $H$ coaction, it is enough that $P$ is an $(H,C)$-bicomodule.

The following lemma is concerned with the existence of  an entwining structure on cotensor product
$A\Box_H P$, induced by the entwining structure on $P$. 
 
  \begin{mlem}\label{entcotlem}
Let $A$, $P$ be, respectively, right and left $H$-comodule algebras, where $H$ is a bialgebra and let 
$C$ be a coalgebra flat as a $\mK$-module. Let $(P,C)_\psi$ be an entwining structure
such that $\psi$ commutes with the left $H$-coaction. Then 
\begin{enumerate}
\item $(A\Box_H P)_{\psi_\Box}$ is an entwining structure
with
\begin{align}
\psi_\Box:C\ot(A\Box_H P)&\rightarrow (A\Box_H P)\ot C, \nonumber\\
c\ot(\sum_ia_i\ot p_i)&\mapsto \sum_i(a_i\ot p_i{}_\alpha)\ot c^\alpha.
\end{align} 
\item If $\psi$ is invertible then $\psi_\Box$ is also invertible with:
\begin{align}
\psi_\Box{}^{-1}:(A\Box_H P)\ot C&\rightarrow C\ot(A\Box_H P), \nonumber\\
(\sum_ia_i\ot p_i)\ot c&\mapsto c_A\ot(\sum_ia_i\ot p_i{}^A).\label{pbidef}
\end{align}
\item If $P\in\Mod^C_P(\psi)$ and $({}^H\coact\ot C)\circ\coact^C(1_P)=1_H\ot\coact^C(1_P)$ then
$A\Box_H P\in \Mod^C_{A\Box_H P}(\psi_\Box)$. 
\item If $P$ is $e$-copointed then $A\Box_H P$ is $e$-copointed.
\end{enumerate} 
\end{mlem}
\begin{proof}
As $C$ is a flat $\mK$-module, 
in order to prove that $\psi_\Box$ has its image in $(A\Box_H P)\ot C$,
it is enough to show that 
$(\coact^H\ot P\ot C)\circ\psi_\Box=(A\ot {}^H\coact\ot C)\circ\psi_\Box$ 
For any
$c\ot (\sum_ia_i\ot p_i)\in C\ot (A\Box_H P)$:
\begin{multline*}
(\coact^H\ot P\ot C)\circ\psi_\Box(c\ot (\sum_ia_i\ot p_i))
=\sum_i(a_i\sw{0}\ot a_i\sw{1}\ot p_i{}_\alpha)\ot c^\alpha\\
=\sum_i(a_i\ot p_i\sw{-1}\ot p_i\sw{0}{}_\alpha)\ot c^\alpha
=\sum_i(a_i\ot p_i{}_\alpha\sw{-1}\ot p_i{}_\alpha\sw{0})\ot c^\alpha\\
=(A\ot {}^H\coact\ot C)\psi_\Box(c\ot (\sum_ia_i\ot p_i)).
\end{multline*}
That $\psi_\Box$ is an entwining follows easily from the definition of $\psi_\Box$. 
Moreover, if  $\psi$ is invertible, $\psi^{-1}$ also commutes with the left $H$-coaction,
and one can prove, with similar computations as those for $\psi_\Box$, that 
the map $\psi_\Box^{-1}$ has its
image in $A\Box_H P$. Obviously $\psi_\Box^{-1}$ is the inverse of $\psi_\Box$.

Assume that $P\in\Mod^C_P(\psi)$
and $({}^H\coact\ot C)\circ\coact^C(1_P)=1_H\ot\coact^C(1_P)$. Then
by Lemma~\ref{hcent}, $P$ is an $(H,C)$-bicomodule, hence $A\Box_H P$ is an algebra and a $C$-comodule (by Lemma \ref{algcot}). Moreover, for any 
$\sum_ia_i\ot p_i,\sum_ja'{}_j\ot p'{}_j\in A\Box_H P$,
\begin{multline*}
\coact^C((\sum_ia_i\ot p_i)(\sum_ja'{}_j\ot p'{}_j))=\sum_{ij}a_ia'{}_j\ot\coact^C(p_ip'{}_j)\\
=\sum_{ij}a_ia'{}_j\ot p_i\sw{0}\psi(p_i\sw{1}\ot p'{}_j)\\
=(\sum_ia_i\ot p_i)\sw{0}\psi_\Box((\sum_ia_i\ot p_i)\sw{1}\ot (\sum_ja'{}_j\ot p'{}_j)),
\end{multline*} 
hence $A\Box_H P\in\Mod^C_{A\Box_H P}(\psi_\Box)$. If $P$ is $e$-copointed, then
\[
\coact^C(1_A\ot 1_P)=1_A\ot\coact^C(1_P)=(1_A\ot 1_P)\ot e=1_{A\Box_H P}\ot e,
\]
hence 
$A\Box_H P$ is $e$-copointed.
\end{proof}

  Recall that a coalgebra $C$ is said to be \emph{coseparable}
   if there exists a $(C,C)$-bicolinear retraction of the comultiplication,
i.e., a map $\delta:C\ot C\rightarrow C$ such that $\delta\circ\comul=C$. If $C$  is a coseparable coalgebra,
then $C$-coalgebra Galois extensions have the following property.
 \begin{mth}(\cite{Brz:Gal}, Theorem 4.6)\label{brz:galsep}
 Let $\mK$ be a field and $P_e(B,C,\psi)$ an $e$-copointed $(P,C)_\psi$-extension with $\psi$ bijective. If $C$ is a coseparable coalgebra and the lifted canonical map
 $\lcan^C_P:P\ot P\rightarrow P\ot C$ is surjective then the canonical map
 $\can^C_P:P\ot_B P\rightarrow P\ot C$ is bijective and $P_e^C(B)$ is a principal $C$-coalgebra Galois extension.
 \end{mth}
 
 \begin{mcor}\label{galsep}
 If $\mK$ is a field, $C$ is a coseparable coalgebra, $P_e(B,C,\psi)$ is an $e$-copointed $(P,C)_\psi$-extension with bijective $\psi$ and there exists a colifting of the translation map
 $\ltrans^C_P:C\rightarrow P\ot P$, then $\can^C_P:P\ot_B P\rightarrow P\ot C$ is bijective and $P_e^C(B)$ is a principal $C$-coalgebra Galois extension.
 \end{mcor}
 \begin{proof} Define map
 \[
 f:P\ot C\ni p\ot c\mapsto p\ltrans^C_P(c)\in P\ot P.
 \]
 By the definition of $\ltrans^C_P$,  $\lcan^C_P\circ f=P\ot C$, hence
 $\lcan^C_P$ must be surjective and conclusion follows from Theorem~\ref{brz:galsep}. 
 \end{proof}

\begin{mlem}\label{equivhco}
 If $P$ is a left $H$-comodule, where $H$ is a Hopf algebra, then for all $\sum_ip_i\ot q_i\in P\ot P$, the following two  conditions are equivalent:
 \begin{gather}
 \sum_i p_i\sw{-1} q_i\sw{-1}\ot p_i\sw{0}\ot q_i\sw{0}=1_H\ot\sum_i p_i\ot q_i,
 \label{hfirst}\\
 \sum_i p_i\sw{-1}\ot p_i\sw{0}\ot q_i=\sum_i Sq_i\sw{-1}\ot p_i\ot q_i\sw{0}.
 \label{hsecond}
 \end{gather}
 \end{mlem}
 \begin{proof}
 
{\noindent\bf (\ref{hfirst})$\Rightarrow$(\ref{hsecond})}
 Apply $(\mul\ot P\ot P)\circ(S\ot{}^H\coact\ot P)$ to both sides of (\ref{hfirst}).
 
{\noindent\bf(\ref{hsecond})$\Rightarrow$(\ref{hfirst})}
 \begin{multline*}
\sum_i p_i\sw{-1}q_i\sw{-1}\ot p_i\sw{0}\ot q_i\sw{0}\\
=\text{[ use (\ref{hsecond}) ] }
\sum_iSq_i\sw{-1}q_i\sw{0}\sw{-1}\ot p_i\ot q_i\sw{0}\sw{0}\\
=\sum_i Sq_i\sw{-1}\sw{1}q_i\sw{-1}\sw{2}\ot p_i\ot q_i\sw{0}
  =1_H\ot\sum_i p_i\ot q_i.
 \end{multline*}
 \end{proof}

 The following two theorems state conditions for the existence, and give an explicit form of the strong connection form.
 
 \begin{mth}\label{thcosepc}
  Assume that:
  \begin{enumerate}
  \item $\mK$ is a field, $H$ is a Hopf algebra, and $C$ is a coseparable coalgebra,
  \item $A$ is  a right $H$-comodule algebra, 
  \item $P$ is an $(H,C)$-bicomodule and a left $H$-comodule algebra. Also, there exists a grouplike 
  $e\in C$ and a bijective entwining $\psi:C\ot P\rightarrow P\ot C$ such that
  $\coact^C(p)=\psi(e\ot p)$ for any $p\in P$,
  \item There exist colifings of the translation maps $\ltrans^H_A:H\rightarrow A\ot A$,
  $\ltrans^C_P:C\rightarrow P\ot P$ which satisfy conditions (\ref{rcovcolift}-\ref{lhcovcol}) 
  Lemma \ref{coliftlem}.
  \end{enumerate}
  Then $(A\Box_H P)^C(R)$, where $R=(A\Box_H P)^{\co C}$, is a principal extension.
 \end{mth}
 \begin{proof} As $\mK$ is a field, any $\mK$-module is $\mK$-flat. 
 Hence all the assumptions of 
 Lemma~\ref{coliftlem} are satisfied, and we know that there exists a colifting of the translation map
 $\ltransc:C\rightarrow (A\Box_H P)\ot (A\Box_H P)$. Moreover the existence
  of $\ltrans^C_P$ which satisfies
 condition (\ref{lhcovcol})  implies, by Lemma~\ref{hcent2} and Lemma~\ref{equivhco} that $\psi$ commutes with the left $H$-coaction.
 Hence, by Lemma~\ref{entcotlem}, there exists an 
 invertible entwining $\psi_\Box$ on $A\Box_H P$ and
 $(A\Box_H P)_e(R,C,\psi_\Box)$ is an $e$-copointed $(A\Box_H P,C)_{\psi_\Box}$-extension. 
 Then the assertion follows by Corollary~\ref{galsep}.
 \end{proof}

 \newcommand{\stlc}{\stl_{A\Box_H P}}
 
 \begin{mth}\label{thstr}
 If 
 \begin{enumerate}
 \item $C$ is a coalgebra and $H$ is a Hopf algebra with a bijective antipode,
 \item $A$ is a right $H$-comodule algebra and $P$ is a left $H$-comodule algebra,
 \item $P_e(B,C,\psi)$ is an $e$-copointed $(P,C)_\psi$-extension with a bijective entwining map, 
 which commutes with the left $H$-coaction,
 \item $C$, $A\ot P$, $A\Box_H P$ are flat  $\mK$-modules,
 \item there exist strong connection forms $\stl_A:H\rightarrow A\ot A$ and
 $\stl_P:C\rightarrow P\ot  P$,
 \item For all $c\in C$,
 \begin{equation}
 c\swa\sw{-1}\ot c\swa\sw{0}\ot c\swb=Sc\swb\sw{-1}\ot c\swa\ot c\swb\sw{0},\label{strcondh}
 \end{equation} 
 where $c\swa\ot c\swb=\stl_P(c)$.
 \end{enumerate}
 Then $(A\Box_H P)_e(R,C,\psi_\Box)$, where $R=(A\Box_H P)^{\co H}$, is an  $e$-copointed\linebreak 
 $(A\Box_H P,C)_{\psi_\Box}$-extension with bijective entwining (Lemma~\ref{entcotlem}), and
 \begin{gather}
\stlc:C\rightarrow (A\Box_H P)\ot (A\Box_H P),\nonumber\\
c\mapsto (c\swb\sw{-1}\swa\ot c\swa)\ot (c\swb\sw{-1}\swb\ot c\swb\sw{0})\label{stlcotfor}
\end{gather}
is a strong connection form.
 \end{mth}
\begin{proof} By Lemma~\ref{hcent}, the  left $H$-coaction commutes with the 
right $C$-coaction. Observe that had we assumed that $P$ is an $(H,C)$-bicomodule, then commuting of  $\psi$ with 
$H$-coaction would follow from (\ref{strcondh}) and Lemmas \ref{equivhco} and \ref{hcent2}.
Hence, by Lemma~\ref{entcotlem},
 $(A\Box_H P)_e(R,C,\psi_\Box)$ is an $e$-copointed  $(A\Box_H P,C)_{\psi_\Box}$-extension
with bijective entwining.  

In particular, $\stl_A$ and $\stl_P$ are coliftings of translation map, which satisfy all of the assumptions
of Lemma~\ref{coliftlem}, hence, by Lemma~\ref{coliftlem}, the 
map $\stlc$ given by  (\ref{stlcotfor}) 
is a well defined colifting of the translation map on $A\Box_H P$. It remains to prove that  (\ref{stlcotfor}) satisfies remaining axioms
(\ref{str1}, \ref{str3}, \ref{str4}) of a strong connection form.
First compute
\begin{multline*}
\stlc(e)=(e\swb\sw{-1}\swa\ot e\swa)\ot (e\swb\sw{-1}\swb\ot e\swb\sw{0})\\
=(1_P\sw{-1}\swa\ot 1_P)\ot (1_P\sw{-1}\swb\ot 1_P\sw{0})=
(1_H\swa\ot 1_P)\ot (1_H\swb\ot 1_P)\\
=(1_A\ot 1_P)\ot (1_A\ot 1_P)=1_{A\Box_H P}\ot 1_{A\Box_H P},
\end{multline*}
Where we use that $\stl_A(e)=1_A\ot 1_A$ and $\stl_P(e)=1_P\ot 1_P$.
Hence, the map $\stlc$ is normalized on $e$ as required for (\ref{str1}).

{\noindent {}}Take any $c\in C$, and compute, 
\begin{multline*}
((A\ot P)\ot\coact^C)\circ\stlc(c)=
(c\swb\sw{-1}\swa\ot c\swa)\ot (c\swb\sw{-1}\swb\ot c\swb\sw{0})\ot c\swb\sw{1}\\
=(c\swb\sw{0}\sw{-1}\swa\ot c\swa)\ot
(c\swb\sw{0}\sw{-1}\swb\ot c\swb\sw{0}\sw{0})\ot c\swb\sw{1}
\\
\text{[by (\ref{str3}) for $\stl_P$] }=
(c\sw{1}\swb\sw{-1}\swa\ot c\sw{1}\swa)\ot (c\sw{1}\swb\sw{-1}\swb\ot c\sw{1}\swb\sw{0})
\ot c\sw{2}\\
=\stlc(c\sw{1})\ot c\sw{2}.
\end{multline*}
Therefore $\stlc$ is a right $C$-comodule map, i.e., the condition (\ref{str3}) is satisfied.

{\noindent{}}Finally, for all $c\in C$,
\begin{multline*}
({}^{C_{\psi_\Box}}\coact\ot (A\ot P))\circ\stlc(c)=
e_A\ot(c\swb\sw{-1}\swa\ot c\swa{}^A)\ot (c\swb\sw{-1}\swb\ot c\swb\sw{0})\\
\text{[use (\ref{str4}) for $\stl_P$] }=
c\sw{1}\ot (c\sw{2}\swb\sw{-1}\swa\ot c\sw{2}\swa)\ot (c\sw{2}\swb\sw{-1}\swb\ot c\sw{2}\swb\sw{0})\\
=c\sw{1}\ot \stlc(c\sw{2}).
\end{multline*}
Therefore, $\stlc$ is a left $C$-comodule map, i.e., the condition (\ref{str4}) is satisfied.
Thus we conclude that $\stlc$ is a strong connection form as required.
\end{proof}

Let $P_e^C(B)$ be a $C$-coalgebra Galois extension, with a strong connection form
$\stl:C\rightarrow P\ot P$, $c\mapsto c\swa\ot c\swb$. Then, obviously, the inverse of the canonical map
is given explicitly by
\[
(\can^C_P)^{-1}(p\ot c)=pc\swa\ot_Bc\swb\ \  \text{for all }c\in C,\ p\in P.
\] 
Hence, explicit formula for the inverse of the canonical map is determined by the explicit formula for 
a strong connection form and by the knowledge of the subalgebra of coinvariants. Next lemma  determines 
a subalgebra of coinvariants of $A\Box_H P$.

\begin{mlem}
Suppose that $H$ is a bialgebra, $C$ is a coalgebra, flat as a $\mK$-module. Let $A$ be a right
$H$-comodule algebra and let $P$ be a left $H$-comodule algebra and an $(H,C)$-bicomodule.
Moreover, suppose that $(P,C)_\psi$ is an entwining structure such that $P$ is an entwined module and 
$\psi$ commutes with the left $H$-coaction. 
Then 
\begin{equation}
(A\Box_H P)^{\co C}=A\Box_H (P^{\co C}).
\end{equation}
\end{mlem}
\begin{proof}
As $P$ is an entwined module, subalgebra of coinvariants is uniquely determined by the exactness
of the following sequence:
\begin{equation}
\xymatrix{
0 \ar[r] & P^{\co C} \ar[r] & P\eqmap{rr}{\coact^C}{p\mapsto p1\sw{0}\ot 1\sw{1}} &  & P\ot C.
}\label{pcois}
\end{equation}
View $P\ot C$ as a left $H$-comodule by left $H$-coaction on the first factor.
The functor $A\Box_H \cdot:{}^H\Mod\rightarrow\Mod_\mK$ is left exact, and, furthermore, 
since $C$ is flat as a $\mK$-module, 
 $A\Box_H(P\ot C)\simeq (A\Box_H P)\ot C$ (c.f. \cite{Tak:GeRed}). Hence, cotensoring sequence
 (\ref{pcois}) with $A$, yields the following exact sequence
 \begin{equation}
 \xymatrix{
0 \ar[r] & A\Box_H (P^{\co C}) \ar[r] & A\Box_H P
\eqmap{rrr}{A\ot \coact^C}{A\ot (p\mapsto p1\sw{0}\ot 1\sw{1})} &  &  &(A\Box_HP)\ot C.
}\label{cotcovseq}
 \end{equation} 
 By Lemmas \ref{algcot} and \ref{entcotlem}, $A\Box_H P$ is an algebra and an entwined module,
 hence sequence (\ref{cotcovseq}) defines uniquely the algebra of coinvariants
 $(A\Box_H P)^{\co C}=A\Box_H (P^{\co C})$.
\end{proof}

\begin{mex}
Let $A$, $P$, $B$, $H$ be as in Example~\ref{matex}, except that we change left $H$-coaction on $P$.
Define the left $\ast$-comodule-algebra $H$-coaction on $P$,
\begin{equation}
{}^H\coact(a)=u^{\ast}\ot a,\ \ {}^H\coact(b)=u\ot b,\ \ {}^H\coact(a^\ast)=u\ot a^\ast,\ \ 
{}^H\coact(b^\ast)=u^{\ast}\ot b^\ast.
\end{equation}
As a $\ast$-algebra, $A\Box_H P$ is generated by,
\begin{align}
\alpha&=a\ot a^\ast,& \beta&=b\ot b, & \gamma&=a\ot b, & \delta&=b\ot a^\ast,\label{mgen}
\end{align}
which satisfy the following commutation relations
\begin{align}
  \alpha\alpha^\ast&=\alpha^\ast\alpha,&
  \beta\beta^\ast&=\beta^\ast\beta,&
  \gamma\gamma^\ast&=\gamma^\ast\gamma,&
  \delta\delta^\ast&=\delta^\ast\delta,\nonumber\\
   \alpha\beta&=\lambda\bar{\lambda'}\beta\alpha,&
  \alpha\beta^\ast&=\bar{\lambda}\lambda'\beta^\ast \alpha,&
  \alpha\gamma&=\bar{\lambda'}\gamma\alpha,&
  \alpha\gamma^\ast&=\lambda'\gamma^\ast\alpha,\nonumber\\
  \alpha\delta&=\lambda\delta\alpha,&
    \alpha\delta^\ast&=\bar{\lambda}\delta^\ast\alpha,&
    \beta\gamma&=\bar{\lambda}\gamma\beta,&
        \beta\gamma^\ast&=\lambda\gamma^\ast\beta,\nonumber\\
        \beta\delta&=\lambda'\delta\beta,&
                \beta\delta^\ast&=\bar{\lambda'}\delta^\ast\beta,&
\gamma\delta&=\lambda\lambda'\delta\gamma,&
\gamma\delta^\ast&=\bar{\lambda}\bar{\lambda'}\delta^\ast\gamma,\label{mcomu}
\end{align}
and, in addition,
\begin{gather}
  \alpha^\ast\alpha+\beta^\ast\beta+\gamma^\ast\gamma+\delta^\ast\delta=1,\label{mrad1}\\
  \alpha\beta=\bar{\lambda'}\gamma\delta,\label{mrad2}
  \end{gather}
  where $\lambda=e^{2\pi i\theta}$ and $\lambda'=e^{2\pi i\theta'}$. 

The algebra of coinvariants
  $R=(A\Box_H P)^{\co H}=A\Box_H B$ is obviously generated by 
  \begin{gather}
  z^1=z\ot 1=\alpha^\ast\alpha+\gamma^\ast\gamma,\ \ 
  z^2=1\ot z=\alpha\alpha^\ast+\delta\delta^\ast,\nonumber\\
  x_+^1=x_+\ot 1=\delta\alpha^\ast+\beta\gamma^\ast,\ \
  x_-^1=x_-\ot 1=\alpha\delta^\ast+\gamma\beta^\ast,\nonumber\\
  x_+^a=aa\ot x_+=\gamma\alpha,\ \ 
  x_-^a=a^\ast a^\ast\ot x_-=\alpha^\ast\gamma^\ast,\nonumber\\
  x_+^b=bb\ot x_+=\beta\delta,\ \ 
  x_-^b=b^\ast b^\ast\ot x_-=\delta^\ast\beta^\ast\nonumber\\
  x_+^{ab}=ab\ot x_+=\lambda'\alpha\beta,\ \
  x_-^{ab}=b^\ast a^\ast\ot x_-=\bar{\lambda'}\beta^\ast\alpha^\ast. \label{mindef}
  \end{gather}
  However  not all of these are 
   independent. Namely,
  \begin{multline}
  x_+^{ab}=\lambda'\alpha\beta=
  \lambda'\alpha\beta(\alpha^\ast\alpha+\beta^\ast\beta+\gamma^\ast\gamma+\delta^\ast\delta)
  \text{ \ [use (\ref{mrad1})]}\\
  =\gamma\delta\alpha^\ast\alpha+\gamma\delta\beta^\ast\beta+
  \lambda'\alpha\beta\gamma^\ast\gamma+\lambda'\alpha\beta\delta^\ast\delta
  \text{\ \ \ [use (\ref{mrad2}) in 1-st \& 2-nd factor]}\\
  =\lambda(\delta\alpha^\ast+\beta\gamma^\ast)\gamma\alpha+
  (\alpha\delta^\ast+\gamma\beta^\ast)\beta\delta
  \text{ \ \ \ [use (\ref{mcomu})]}\\
  =\lambda x_+^1x_+^a+x_-^1x_+^b 
  \text{\ \ \ [use (\ref{mindef})]},
  \end{multline}
  and
  \begin{equation}
  x_-^{ab}=(x_+^{ab})^\ast=
  \bar{\lambda}x_-^ax_-^1+x_-^bx_+^1.
  \end{equation}
  The remaining generators satisfy the following commutation relations:
  \begin{align}
  x_+^1x_-^1&=x_-^1x_+^1, & x_+^ax_-^a&=x_-^ax_+^a, & x_+^bx_-^b&=x_-^bx_+^b,\nonumber\\
  x_+^1x_+^a&=\bar{\lambda}^2x_+^ax_+^1, & x_+^1x_-^a&=\lambda^2x_-^ax_+^1, &
  x_+^1x_+^b&=\bar{\lambda}^2x_+^bx_+^1,\nonumber\\
  x_+^1x_-^b&=\lambda^2x_-^bx_+^1, & x_+^ax_+^b&=\lambda^2x_+^bx_+^a, &
  x_+^ax_-^b&=\bar{\lambda}^2x_-^bx_+^a,\label{stcomu}
  \end{align}
  and $z^1$ and $z^2$ are central in $R$. In addition to
  relations (\ref{stcomu}), the generators $z^1$, $z^2$, $x_+^1$, $x_-^1$,
  $x_+^a$, $x_-^a$, $x_+^b$, $x_-^b$ satisfy
  \begin{gather}
  x_+^1x_-^1+(z^1)^2=z^1,\label{eq1}\\ 
  x_+^ax_-^a=(z^1)^2z^2(1-z^2),\label{eq2}\\
  x_+^bx_-^b=(1-z^1)^2z^2(1-z^2),\label{eq3} \\
   x_+^ax_-^b=\bar{\lambda}(x_-^1)^2z^2(1-z^2).\label{eq4}
  \end{gather}
  
  To gain better geometric understanding of the algebra $A\Box_H B$ we look at the 
  classical case, whereby  $\theta=\theta'=0$.
  In accordance with remarks in the introduction about the classical interpretation of the cotensor product,
  the algebra of coinvariants $A\Box_H B$ is a commutative  algebra of functions on
  a fibered space with the base 
  $S^2$ and the fibre $S^2$. Denote this space by $S^2_{S^2}$, i.e., 
  \[
  S^2_{S^2}=\{(x_+^1,x_+^a,x_+^b,z^1,z^2)\in\mC\times\mC\times\mC\times\mR\times\mR|
  \text{ 
  eq's (\ref{eq1})-(\ref{eq4})
  are satisfied}\},
  \]
  where we set $\lambda=\lambda'=1$ in equations  (\ref{eq1})-(\ref{eq4}).
  It is easy to see that parameters $x^1_{\pm}$, $z^1$ and (\ref{eq1}) describe
  the base space of this fibration. Hence we can define a surjection on the base space
  $\pi:S^2_{S^2}\rightarrow S^2$, by 
  \[
  \pi:(x_+^1,x_+^a,x_+^b,z^1,z^2)\mapsto (x_+^1,z).
  \]
  Define also $U_0,U_1\subset S^2$, $U_k=\{(x_+,z)\in S^2|z\neq k\}$.
  
  For any parameter $\alpha\in\mR$, $\alpha>0$, the equation
  \[
  x_+x_-=\alpha z(1-z),
  \]
  where $z\in\mR$, $x_-=\overline{x_+}\in\mC$,
 describes an ellipsoid with equatorial radius $\frac{\sqrt{\alpha}}{2}$, and longitudal diameter $1$,
 which is obviously homeomorphic with $S^2$ by the substitution 
 $(x_+,z)\mapsto (\frac{x_+}{\sqrt{\alpha}},z)$. Define two continuous maps, by 
 $\phi_0:\pi^{-1}(U_0)\rightarrow U_0\times S^2$ and
  $\phi_1:\pi^{-1}(U_1)\rightarrow U_0\times S^2$,
\begin{gather*}
\phi_0:(x_+^1,x_+^a,x_+^b,z^1,z^2)\mapsto ((x_+^1,z^1),(\frac{x_+^a}{z^1},z^2)),\\
\phi_1:(x_+^1,x_+^a,x_+^b,z^1,z^2)\mapsto ((x_+^1,z^1),(\frac{x_+^b}{1-z^1},z^2)).
\end{gather*}
These maps have the inverses 
\begin{equation*}
(\phi_0)^{-1}:((x_+^1,z^1),(x_+,z))\mapsto (x_+^1,x_+^a,x_+^b,z^1,z^2)=
(x_+^1,z^1x_+,(x_+^1)^2\frac{x_+}{z^1},z^1,z),
\end{equation*}
and
\begin{equation*}
\begin{split}
(\phi_1)^{-1}:((x_+^1,z^1),(x_+,z))&\mapsto (x_+^1,x_+^a,x_+^b,z^1,z^2)\\
&\quad{}=(x_+^1,(x_-^1)^2\frac{x_+}{1-z^1},(1-z^1)x_+,z^1,z),
\end{split}
\end{equation*}
which are also continuous. Hence $\phi_0$ and $\phi_1$ are homeomorphisms which define 
local trivialisations of the fibre bundle $S^2_{S^2}$. Moreover they cannot be extended to the whole
of $S^2_{S^2}$, hence $S^2_{S^2}$ is a nontrivial fibration.
 $A\Box_H B$ is a quantum deformation of algebra 
of functions on a nontrivial fibre bundle $S^2_{S^2}$.

 For any $\theta\in\mR$. $C^0_\theta(S^3)$ admits a strong connection form given explicitly by
 \begin{gather}
 \stl(u^n)=
 \sum_{m=0}^n {n\choose m} b^\ast{}^ma^\ast{}^{n-m}\ot a^{n-m} b^m,\label{ms1}\\
 \stl(u^\ast{}^n)=
  \sum_{m=0}^n {n \choose m}b^ma^{n-m}\ot a^\ast{}^{n-m} b^\ast{}^m.\label{ms2}
 \end{gather}
for all $n\in\mN$ (c.f. \cite{BrzSit:Mat}). Hence $A=C^0_\theta(S^3)$ and $P=C^0_{\theta'}(S^3)$ admit strong connection forms $\stl_A$ and $\stl_P$, given by the above formulae. Observe that
$\stl_P$ satisfies condition (\ref{strcondh}) of Theorem~\ref{thstr}. Hence, by 
Theorem~\ref{thstr}, there exists a strong connection form $\stl_{A\Box_H P}$ on $A\Box_H P$. 
Explicitly, by (\ref{stlcotfor}), (\ref{ms1}) and (\ref{ms2}), for any $n\in\mN$,
\begin{equation}
\begin{split}
\stlc(u^n)&=\sum_{m=0}^{\lfloor\frac{n}{2}\rfloor}\sum_{k=0}^{n-2m}
{n\choose m}{n-2m\choose k}(b^ka^{n-2m-k}\ot b^\ast{}^ma^\ast{}^{n-m})\\
&\quad\quad{}\ot (a^\ast{}^{n-2m-k}b^\ast{}^k\ot a^{n-m} b^m)\\
&\quad{}+\sum_{m=\lfloor\frac{n}{2}\rfloor+1}^n\sum_{k=0}^{2m-n}
{n \choose m}{2m-n\choose k}
(b^\ast{}^ka^\ast{}^{2m-n-k}\ot b^\ast{}^ma^\ast{}^{n-m})\\
&\quad\quad{}\ot (a^{2m-n-k}b^k\ot a^{n-m}b^m),
\end{split}\label{strmaexpl1}
\end{equation}
and
\begin{equation}
\begin{split}
\stlc(u^{-n})&=
\sum_{m=0}^{\lfloor\frac{n}{2}\rfloor}\sum_{k=0}^{n-2m}
{n\choose m}{n-2m\choose k}
(b^\ast{}^ka^\ast{}^{n-2m-k}\ot b^ma^{n-m})\\
&\quad\quad {}\ot (a^{n-2m-k}b^k\ot a^\ast{}^{n-m}b^\ast{}^m)\\
&\quad{}+\sum_{m=\lfloor\frac{n}{2}\rfloor+1}^n\sum_{k=0}^{2m-n}
{n \choose m}{2m-n\choose k}
(b^ka^{2m-n-k}\ot b^ma^{n-m})\\
&\quad\quad{}\ot (a^\ast{}^{2m-n-k}b^\ast{}^k\ot a^\ast{}^{n-m}b^\ast{}^m).
\end{split}
\label{strmaexpl2}
\end{equation}
In order to express the above formulae in terms of generators (\ref{mgen}), observe that
\begin{multline*}
 a^\ast{}^{n-2m-k}b^\ast{}^k\ot a^{n-m} b^m=a^\ast{}^{n-2m-k}b^\ast{}^k(a^\ast a+b^\ast b)^m
 \ot a^{n-m} b^m\\
 =\sum_{t=0}^m{m\choose t} \bar{\lambda}^{t(k+m-t)}
 a^\ast{}^{n-2m-k+t} b^\ast{}^{m+k-t}a^tb^{m-t}\ot a^{n-m}b^m\\
 =\sum_{t=0}^m{m\choose t} \bar{\lambda}^{t(k+m-t)}
 \alpha^\ast{}^{n-2m-k+t}\delta^\ast{}^{m+k-t}\gamma^{t}\beta^{m-t}.
 \end{multline*}
 Similarly, consider the second factor of the tensor product in the second summand in (\ref{strmaexpl1}).
 If $0\leq k\leq 2m-n$ and $\lfloor\frac{n}{2}\rfloor<m\leq n$, then $k\leq m$. It follows, that
 \begin{multline*}
 a^{2m-n-k}b^k\ot a^{n-m}b^m=
 (a^\ast a+b^\ast b)^{n-m}a^{2m-n-k}b^k\ot a^{n-m}b^{m}\\
 =\sum_{t=0}^{n-m}{n-m\choose t}a^\ast{}^tb^\ast{}^{n-m-t}b^{n-m-t}
 a^{2m-n-k+t}b^k\ot a^{n-m}b^{m}\\
 =\sum_{t=0}^{n-m}{n-m\choose t}\lambda^{(2m-n-k+t)k}
 \alpha^\ast{}^t\delta^\ast{}^{n-m-t}\beta^{n-m-t+k}\gamma^{2m-n-k+t}.
 \end{multline*}
 Hence
 \begin{equation}
 \begin{split}
 \stlc(u^n)&=\sum_{m=0}^{\lfloor\frac{n}{2}\rfloor}\sum_{k=0}^{n-2m}\sum_{t=0}^m\sum_{s=0}^m
{n\choose m}{n-2m\choose k}{m\choose t}{m\choose s}\lambda^{(k+m)(t-s)-t^2+s^2}\\
&\quad\quad{}\cdot\beta^\ast{}^{m-t}\gamma^\ast{}^{t}\delta^{k+m-t}\alpha^{n-2m-k+t}
\ot\alpha^\ast{}^{n-2m-k+s}\delta^\ast{}^{k+m-s}\gamma^s\beta^{m-s}\\
&\quad{}+\sum_{m=\lfloor\frac{n}{2}\rfloor+1}^{n}\sum_{k=0}^{2m-n}\sum_{t=0}^m\sum_{s=0}^m
{n\choose m}{2m-n\choose k}{m\choose t}{m\choose s}\bar{\lambda}^{k(t-s)}\\
&\quad\quad{}\cdot\gamma^\ast{}^{2m-n-k+t}\beta^\ast{}^{n-m+k-t}\delta^{n-m-t}\alpha^t\\
&\quad\quad\quad{}\ot \alpha^\ast{}^s\delta^\ast{}^{n-m-s}\beta^{n-m-s+k}\gamma^{2m-n-k+s}.
 \end{split}
 \end{equation}
 Similarly we prove that 
 \begin{equation}
 \begin{split}
 \stlc(u^{-n})&=\sum_{m=0}^{\lfloor\frac{n}{2}\rfloor}\sum_{k=0}^{n-2m}\sum_{t=0}^m\sum_{s=0}^m
{n\choose m}{n-2m\choose k}{m\choose t}{m\choose s}\lambda^{(k+m)(t-s)-t^2+s^2}\\
&\quad\quad{}\cdot
\alpha^\ast{}^{n-2m-k+s}\delta^\ast{}^{k+m-s}\gamma^s\beta^{m-s}\ot
\beta^\ast{}^{m-t}\gamma^\ast{}^{t}\delta^{k+m-t}\alpha^{n-2m-k+t}\\
&\quad{}+\sum_{m=\lfloor\frac{n}{2}\rfloor+1}^{n}\sum_{k=0}^{2m-n}\sum_{t=0}^m\sum_{s=0}^m
{n\choose m}{2m-n\choose k}{m\choose t}{m\choose s}\bar{\lambda}^{k(t-s)}\\
&\quad\quad{}\cdot
\alpha^\ast{}^s\delta^\ast{}^{n-m-s}\beta^{n-m-s+k}\gamma^{2m-n-k+s}\\
&\quad\quad\quad{}\ot
\gamma^\ast{}^{2m-n-k+t}\beta^\ast{}^{n-m+k-t}\delta^{n-m-t}\alpha^t.
 \end{split}
 \end{equation}
 As $H=C^0(U(1))$ is a coseparable Hopf algebra, by Theorem~\ref{thcosepc},\linebreak 
 $(A\Box_H P)^H(A\Box_H B)$ is an $H$-Hopf Galois extension.
\end{mex}

\section*{Acknowledgements}
I would like to thank Tomasz Brzezi\'nski for many fruitful discussions. My research of 
is supported by the EPSRC grant~GR/S01078/01.

\bibliographystyle{plain}

 \end{document}